\documentclass[preprint,11pt,a4paper]{elsarticle}
\usepackage[T1]{fontenc}
\usepackage{lmodern}
\usepackage{geometry}
\usepackage{amsmath,amssymb,amsthm,mathtools,enumitem,booktabs,array,tabularx}
\usepackage[colorlinks=true,linkcolor=blue,urlcolor=blue,citecolor=blue]{hyperref}
\usepackage{microtype}

\usepackage{indentfirst}
\setlist{itemsep=3pt,topsep=5pt}
\let\paragraph\elsparagraph
\newtheorem{theorem}{Theorem}[section]
\newtheorem{lemma}[theorem]{Lemma}
\newtheorem{proposition}[theorem]{Proposition}
\newtheorem{corollary}[theorem]{Corollary}
\theoremstyle{definition}
\newtheorem{definition}[theorem]{Definition}
\newtheorem{remark}[theorem]{Remark}

\newcommand{\G}{\mathrm{GpsUL}^{*}}
\newcommand{\Gtwo}{\mathrm{GpsUL}^{**}}

\newcommand{\seq}{\Rightarrow}
\newcommand{\Om}{\mathsf A_{\mathrm{ord}}}
\newcommand{\ECs}{\mathrm{EC}_{\mathrm{coh}}}
\newcommand{\Par}{\operatorname{Par}}
\newcommand{\supp}{\operatorname{supp}}
\newcommand{\lca}{\operatorname{lca}}
\newcommand{\Cut}{\mathrm{CUT}}
\newcommand{\WCT}{\mathrm{WCT}}
\newcommand{\WCM}{\mathrm{WCM}}
\newcommand{\COM}{\mathrm{COM}}
\newcommand{\EC}{\mathrm{EC}}
\newcommand{\EW}{\mathrm{EW}}
\journal{Fuzzy Sets and Systems}
\biboptions{numbers,sort&compress}
\date{}
\hypersetup{pdftitle={Density elimination for weakly commutative pseudo-uninorm logic:
a correction and an explicit theory-level completeness proof},pdfauthor={S.-M. Wang}}
\begin{document}
\begin{frontmatter}
\title{Density elimination for weakly commutative pseudo-uninorm logic:\newline
a correction and an explicit theory-level completeness proof}
\author{S.-M. Wang\corref{cor1}}
\ead{wangsanmin@hotmail.com}
\affiliation{organization={Faculty of Science, Zhejiang Sci-Tech University},
             city={Hangzhou},
             postcode={310018},
             country={China}}
\cortext[cor1]{Corresponding author. Tel.: +86-136-5581-2543.}
\begin{abstract}
We correct the density-elimination argument for weakly commutative
pseudo-uninorm logic in Wang, \emph{Symmetry} 11 (2019), 368,
retaining its basic ideas and overall strategy.
The published definition of a labelled copy permits independent
renamings of antecedent and succedent labels.
Under the syntactic splicing interpretation used in that argument,
a contraction-free derivable hypersequent and a four-element
residuated chain refute the claimed compatibility of the resulting
contraction with ordered density elimination.
We require a single renaming on both sides of each closed block
and prove that the linear rules and coherent closed-block
contraction preserve provability of every member of the associated
Cartesian families.
An occurrence-based separation argument, together with normalization
to cuts on the monoidal unit, establishes strong density elimination.
We then derive standard completeness from density, retaining the
original conclusion for formula theorems.
For an explicit Hilbert presentation with independent conjugation
parameters, we give a direct relative chain-countermodel proof
for arbitrary theories in a countable language.
Combined with the full bounded embedding of Wang (2026), this
recovers the corresponding strong standard completeness result
already stated there.
The distinguished constant $f$, when present, is interpreted
independently of the monoidal unit and the lattice bounds.
The main results of this paper have been formally verified in Lean~4.
\end{abstract}
\end{frontmatter}

\section{Introduction}

Noncommutative fuzzy logics based on residuated chains require care
in passing from algebraic completeness to standard completeness.
Wang and Zhao~\cite{WZ2009} showed that $\mathrm{HpsUL}$ is not
complete with respect to pseudo-uninorm algebras on the real unit
interval. To address this standard-completeness problem,
Wang~\cite{Wang2013} introduced the weakly commutative extension
$\mathrm{HpsUL}^{*}$.

Wang~\cite{Wang2019} subsequently proposed a proof of standard completeness
for $\mathrm{HpsUL}^{*}$ through density elimination in its hypersequent
calculus $\G$. This follows the general syntactic strategy of Metcalfe
and Montagna~\cite{MM2007}: prove completeness after adjoining density,
then eliminate density. The noncommutative proof uses ordered splitting
and splicing and a separation construction related to the earlier
labelled method in~\cite{Wang2015}. The wider algebraic setting of density
is discussed in~\cite{MT2017}.

The obstruction examined here is internal to the labelled proof.
Definition~5(iv) of~\cite{Wang2019} allows the antecedent and succedent
occurrences of one label to be renamed independently. The resulting
relation does not preserve the producer--consumer incidence that determines
an ordered splicing result. Section~\ref{def:section} gives a concrete
counterexample to the contraction--density compatibility asserted in
Lemma~5 when the splicing operation is applied to the contraction
conclusion as in the published proof. The example satisfies the
provable-premise hypothesis and has a unique terminal splicing output. It is not a counterexample to the
standard-completeness statement.

The present density-elimination argument follows the basic ideas
and overall strategy of the original paper~\cite{Wang2019};
the corrections concern technical aspects of the proof and supply
the missing justifications needed to carry out that strategy rigorously.
The local invariant requires provability of every member of every
actual expansion. By weak exchange, it suffices to prove all members of the
Cartesian family associated with one fixed choice of cycle cuts. Coherent contraction then uses a
diagonal choice of corresponding endpoints and ordinary external contraction. The global construction tracks actual
contraction histories and separates sources by a well-founded induction
on antichain size and occurrence counts. Neither a prescribed intermediate
contraction schedule nor uniqueness of separated fragments is needed.

Theorems~\ref{thm:strong-density} and~\ref{thm:strong-standard}, together
with Corollary~\ref{cor:fss-weak-standard}, are the main conclusions.
The density argument retains the density-elimination and weak
standard-completeness conclusions of~\cite{Wang2019}. To avoid
an unproved normalization premise, Appendix~\ref{uc:section} establishes
the reduction of all cuts to cuts on the unit. Appendices
\ref{sec:ordered-preparation} and~\ref{sec:ordered-separation} give the
complete preparation and separation constructions. The dense-countermodel and completion arguments in
Section~\ref{sec:fss-semantics} establish completeness of the
density-augmented calculus. The passage back to the original calculus
uses density elimination. All auxiliary syntactic results needed for the density repair
are proved here; the original constrained-contraction lemma and its
separation algorithm are not used as black boxes.

A corresponding strong standard completeness result for arbitrary
theories already appears in~\cite[Corollary~10.7]{Wang2026}.
Section~\ref{sec:strong-standard} supplies a direct relative
chain-countermodel proof for arbitrary theories in a countable
language of the independent-parameter Hilbert calculus specified
in Appendix~\ref{sec:hilbert-explicit}. This proof fixes the
consequence relation explicitly, constructs a maximal theory avoiding a
proposed conclusion, and shows that its relative Lindenbaum algebra is a
chain. It then combines this logical construction with the full bounded
standard embedding of~\cite[Corollary~10.5]{Wang2026}. The contribution here
is the explicit relative proof for this specified presentation, together
with the density repair; the corresponding strong completeness conclusion
is not claimed as a first result. The full representation theorem is an
external algebraic input and is not reproved in this manuscript.

Finite syntactic support is proved before the maximal-theory construction;
semantic compactness follows from strong completeness and is not assumed.
The original $\G$--Hilbert correspondence remains a correspondence for
formula theorems. Table~\ref{tab:proof-dependencies} separates the two
completeness routes and their inputs.

\begin{table}[htbp]
\centering\small
\caption{Proof dependencies.}
\label{tab:proof-dependencies}
\begin{tabularx}{\textwidth}{@{}>{\raggedright\arraybackslash}p{0.24\textwidth}>{\raggedright\arraybackslash}X >{\raggedright\arraybackslash}p{0.28\textwidth}@{}}
\toprule
Conclusion & Required input & Location and status\\
\midrule
Strong density & Unit-cut normalization, preparation, separation and local
Cartesian compatibility & Proved here: Appendices~\ref{uc:section}--\ref{sec:ordered-separation}
and Sections~\ref{sec:cartesian-density}--\ref{sec:strong-density}\\[3pt]
Weak standard completeness & Density elimination and completion of the
constructed dense countermodel & Proved here: Section~\ref{sec:fss-semantics}\\[3pt]
Relative chain countermodel & Finite proof support, relative algebra and
comparison extension & Proved here: Section~\ref{sec:strong-standard}\\[3pt]
Full bounded standard embedding &
Representation of nontrivial countable bounded unit-cyclic
residuated chains, without a density assumption &
External input:
\cite[Corollary~10.5]{Wang2026}\\[3pt]
Strong standard completeness & Relative chain countermodel and full bounded
standard embedding & Derived here for the explicit calculus; corresponding
result in~\cite[Corollary~10.7]{Wang2026}\\
\bottomrule
\end{tabularx}
\end{table}

The main results of this paper have been formally verified in Lean~4.
The accompanying source code and reproducibility materials are available
in the archived release \emph{GpsUL-1.0.0}~\cite{Wang2026GpsUL}.
The release contains the Lean source code, the corresponding manuscript,
build and axiom-audit scripts, and a statement-level correspondence map.
Section~\ref{sec:scope} describes the proof-development and verification
methods.

\section{Calculi and notation}\label{sec:calculus}

Formulas are built from countably many variables, the constants
$t,\bot,\top$, and $\cdot,\backslash,/,\wedge,\vee$. We use the residual
notation $U\backslash V$ and $V/U$; these are $U\to V$ and
$U\rightsquigarrow V$, respectively, in~\cite{Wang2019}. Both notations
will be used when describing an original rule. Antecedents are finite
ordered words, including the empty word. Hypersequents are finite nonempty
multisets of single-conclusion sequents. Vertical bars denote external
multiset union; an external context can be empty. Equal formulas or
components can have distinct occurrence identities. The formal syntax
represents the external multiset by a list with a rule for permuting
component occurrences. Its empty list serves as a context, and induction
on an ordinary proof gives a nonempty conclusion. For the saturation
construction in Section~\ref{sec:fss-semantics} we use a countably
infinite supply of variables. Section~\ref{sec:strong-standard} treats
any countable alphabet, including finite alphabets.

A bounded residuated chain has an associative product with unit $e$ and
residuals satisfying
\[
 xy\le z\quad\Longleftrightarrow\quad y\le x\backslash z
          \quad\Longleftrightarrow\quad x\le z/y.
\]
A sequent holds under a valuation when the ordered product of its
antecedent, with empty product $e$, is below its succedent. A hypersequent
holds under a fixed valuation when at least one component holds under
that valuation. Validity in an algebra means that this condition holds
for every valuation. We interpret $t$ by $e$. Thus validity of a formula
$F$ means $e\le v(F)$ for every valuation $v$. A standard algebra has lattice reduct
$[0,1]$, and no strict-interiority condition on $e$ is added.
The original hypersequent result concerns formula theorems;
Section~\ref{sec:strong-standard} defines and proves strong standard
completeness for Hilbert consequence from arbitrary theories.

We use the hypersequent presentation of~\cite{MOG2009,Wang2019}.
The initial sequents of $\G$ are
\[
 A\seq A,\qquad \seq t,\qquad
 \Gamma,\bot,\Delta\seq A,\qquad\Gamma\seq\top.
\]
The external rules are
\[
 \frac{G\mid S\mid S}{G\mid S}\;\EC,
 \qquad \frac{G}{G\mid S}\;\EW,
\]
\begin{equation}\label{eq:cart-com-rule}
 \frac{G_1\mid\Gamma_1,\Pi_1,\Delta_1\seq A_1\qquad
       G_2\mid\Gamma_2,\Pi_2,\Delta_2\seq A_2}
 {G_1\mid G_2\mid\Gamma_1,\Pi_2,\Delta_1\seq A_1
                  \mid\Gamma_2,\Pi_1,\Delta_2\seq A_2}\;\COM.
\end{equation}
The logical rules are
\begin{gather*}
 \frac{G\mid\Gamma,U,V,\Delta\seq B}
 {G\mid\Gamma,U\cdot V,\Delta\seq B}\;(\cdot L),\qquad
 \frac{G_1\mid\Pi\seq U\qquad G_2\mid\Sigma\seq V}
 {G_1\mid G_2\mid\Pi,\Sigma\seq U\cdot V}\;(\cdot R),\\[3pt]
 \frac{G_1\mid\Gamma,V,\Delta\seq B\qquad G_2\mid\Sigma\seq U}
 {G_1\mid G_2\mid\Gamma,\Sigma,U\backslash V,\Delta\seq B}
 \;(\backslash L),\qquad
 \frac{G\mid U,\Pi\seq V}{G\mid\Pi\seq U\backslash V}
 \;(\backslash R),\\[3pt]
 \frac{G_1\mid\Sigma\seq U\qquad G_2\mid\Gamma,V,\Delta\seq B}
 {G_1\mid G_2\mid\Gamma,V/U,\Sigma,\Delta\seq B}\;(/L),\qquad
 \frac{G\mid\Pi,U\seq V}{G\mid\Pi\seq V/U}\;(/R),\\[3pt]
 \frac{G_1\mid\Gamma\seq U\qquad G_2\mid\Gamma\seq V}
 {G_1\mid G_2\mid\Gamma\seq U\wedge V}\;(\wedge R),\\[3pt]
 \frac{G\mid\Gamma,U,\Delta\seq B}
 {G\mid\Gamma,U\wedge V,\Delta\seq B}\;(\wedge L_1),\qquad
 \frac{G\mid\Gamma,V,\Delta\seq B}
 {G\mid\Gamma,U\wedge V,\Delta\seq B}\;(\wedge L_2),\\[3pt]
 \frac{G_1\mid\Gamma,U,\Delta\seq B\qquad G_2\mid\Gamma,V,\Delta\seq B}
 {G_1\mid G_2\mid\Gamma,U\vee V,\Delta\seq B}\;(\vee L),\\[3pt]
 \frac{G\mid\Gamma\seq U}{G\mid\Gamma\seq U\vee V}\;(\vee R_1),\qquad
 \frac{G\mid\Gamma\seq V}{G\mid\Gamma\seq U\vee V}\;(\vee R_2),\qquad
 \frac{G\mid\Gamma,\Delta\seq B}{G\mid\Gamma,t,\Delta\seq B}\;(tL).
\end{gather*}
Finally, $\G$ has weak commutativity and unrestricted cut:
\[
 \frac{G\mid\Gamma,\Delta\seq t}{G\mid\Delta,\Gamma\seq t}\;\WCM,
\qquad
 \frac{G_1\mid\Gamma,A,\Delta\seq B\qquad G_2\mid\Pi\seq A}
 {G_1\mid G_2\mid\Gamma,\Pi,\Delta\seq B}\;\Cut.
\]
The calculus $\Gtwo$ replaces unrestricted cut by its instance with
$A=t$, denoted WCT. For reference, the two unit-specific rules are
\begin{equation}\label{eq:cart-wcm-wct}
 \frac{G\mid\Gamma,\Delta\seq t}{G\mid\Delta,\Gamma\seq t}\;\WCM,
\qquad
 \frac{G_1\mid\Gamma,t,\Delta\seq B\qquad G_2\mid\Pi\seq t}
 {G_1\mid G_2\mid\Gamma,\Pi,\Delta\seq B}\;\WCT.
\end{equation}
The terminology refers to the same unit condition in three presentations.
The rule WCM is called \emph{weak exchange}; the Hilbert axiom
$(A\rightsquigarrow t)\to(A\to t)$ is called \emph{weak commutativity}.
Its algebraic form $e/a\le a\backslash e$ is equivalent, by residuation
and interchanging the two variables, to \emph{unit cyclicity}
\[
                 ab\le e\quad\Longleftrightarrow\quad ba\le e.
\]
None of these conditions asserts $ab=ba$.
No internal structural rule is included. In particular, preserving the
positions of inserted words is essential in every construction below.

Write $\G_D$ for $\G$ extended by the density rule
\begin{equation}\label{eq:main-density}
 \frac{G\mid\Pi\seq p\mid\Gamma,p,\Delta\seq B}
 {G\mid\Gamma,\Pi,\Delta\seq B}\;(D),
 \qquad p\notin\operatorname{Var}(G,\Gamma,\Pi,\Delta,B).
\end{equation}
The characteristic variable $p$ is fresh for all the indicated data.
Later labels $p_i$ track occurrences of this distinguished variable; they
are not a family of independent semantic variables.

\section{The copy relation and the contraction--density defect}
\label{def:section}

All original numbering in this section refers to the published thirteen-page
version of Wang~\cite{Wang2019}. Definition~5(iv), on pp.~4--5, requires
the two blocks to have disjoint label supports. A copy is witnessed by
one bijection between their antecedent label sets and another between
their succedent label sets, together with a matching of component
occurrences that realizes these replacements at the displayed positions.
This convention
preserves the unlabelled sequents but can change which producer is attached
to which antecedent position.  In a noncommutative calculus that attachment
determines the ordered result of splicing.  The following counterexample tests the universal claim in Lemma~5
using the syntactic splicing clause on its contracted conclusion.

Write $\mathrm{GL}^{\Omega}_{\mathrm{cf}}$ for the calculus of Section~3 of
\cite{Wang2019}, without its constrained external contraction.  Put
\[
 C=a\cdot b,\qquad D=b\cdot a,\qquad E=C\wedge D,
\]
where $a,b$ are ordinary propositional variables.  For an eigenvariable-free
formula $X$, let
\begin{align}
 U_X&=(p_1,p_2\Rightarrow X)\mid(a\Rightarrow p_1)
                           \mid(b\Rightarrow p_2),\label{def:U}\\
 V_X&=(p_3,p_4\Rightarrow X)\mid(a\Rightarrow p_4)
                           \mid(b\Rightarrow p_3).\label{def:V}
\end{align}
The subscripts on $p$ are identification numbers of occurrences; they are
not independent propositional variables of the underlying unlabelled
semantics.

\begin{proposition}\label{def:counter}
There is a $\mathrm{GL}^{\Omega}_{\mathrm{cf}}$-provable closed hypersequent
$H$ to which Construction~3(i) of~\cite{Wang2019} applies. On extending
the ordered splicing clause of Definition~9(ii) syntactically to its
contracted conclusion, the unique terminal output is unprovable in
$\mathrm{GpsUL}^{\Omega}$. Thus the contraction--density compatibility
asserted in the original Lemma~5 fails on this reading of its density
image.
\end{proposition}

\begin{proof}
First, $U_C$ and $V_D$ are provable without any contraction.  The two
identities for $a,b$ and $\cdot_r$ give $a,b\Rightarrow C$.  Communication
with fresh labelled identities gives
\[
 \frac{a,b\Rightarrow C\qquad p_1\Rightarrow p_1}
 {p_1,b\Rightarrow C\mid a\Rightarrow p_1}\; (\mathrm{COM}),
\]
\[
 \frac{p_1,b\Rightarrow C\mid a\Rightarrow p_1
                   \qquad p_2\Rightarrow p_2}
 {p_1,p_2\Rightarrow C\mid a\Rightarrow p_1\mid b\Rightarrow p_2}
 \;(\mathrm{COM}).
\]
Similarly, starting with $b,a\Rightarrow D$, communication gives
\[
 \frac{b,a\Rightarrow D\qquad p_3\Rightarrow p_3}
 {p_3,a\Rightarrow D\mid b\Rightarrow p_3}\; (\mathrm{COM}),
\]
\[
 \frac{p_3,a\Rightarrow D\mid b\Rightarrow p_3
                   \qquad p_4\Rightarrow p_4}
 {p_3,p_4\Rightarrow D\mid b\Rightarrow p_3\mid a\Rightarrow p_4}
 \;(\mathrm{COM}).
\]
Every two-premise inference has closed, label-disjoint premises.  Applying
$\wedge_{rw}$, as defined in Definition~5(vii) of \cite{Wang2019}, to the
two components with ordinary succedents yields
\begin{equation}\label{def:meet}
 \frac{U_C\qquad V_D}{H:=U_E\mid V_E}\;(\wedge_{rw}).
\end{equation}
This is a proof in $\mathrm{GL}^{\Omega}_{\mathrm{cf}}$; it uses neither
weakening, contraction, cut nor weak exchange.

The left and right label sets of $U_E$ are both $\{1,2\}$, and those of
$V_E$ are both $\{3,4\}$. Thus the blocks are closed and disjoint.
In either block a closed subhypersequent containing a producer must
contain its consumer, the root. A closed subhypersequent containing
the root must contain both producers. Consequently each block is a
minimal nonempty closed unit. Each has a root with ordinary succedent
$E$ and two labelled producers, so both are splicing units in the
sense of Definition~10 of~\cite{Wang2019}. The two bijections in
Definition~5(iv) are witnessed by
\begin{equation}\label{def:maps}
 \sigma_l(1)=3,\quad\sigma_l(2)=4,
 \qquad\sigma_r(1)=4,\quad\sigma_r(2)=3.
\end{equation}
The matching takes root to root and matches producers by their
antecedents $a$ and $b$. In Construction~3 the outside hypersequent is
$G'=\varnothing$. Thus there is no outside splitting unit for which
the two copies would have to lie on the same side of a split. Case~(i)
applies to these two splicing units and takes $G''=G'$. Its precise
instance is
\begin{equation}\label{def:contract}
 \frac{U_E\mid V_E}{U_E}\;\langle\mathrm{EC}^{*}_{\Omega}\rangle.
\end{equation}
In particular, the provable-premise hypothesis of Lemma~5 is met.

For clarity, write $\mathcal D_{\mathrm{sp}}(U_E)$ for the set of terminal
outputs obtained by applying the word replacement in Definition~9(ii)
to the remaining labels of $U_E$. Each producer has a unique consumer
position in the ordinary root. Either label may be spliced first. The
remaining label is then forced, giving the two paths
\[
\begin{aligned}
 U_E&\xrightarrow{1} (a,p_2\Rightarrow E)\mid(b\Rightarrow p_2)
       \xrightarrow{2} (a,b\Rightarrow E),\\
 U_E&\xrightarrow{2} (p_1,b\Rightarrow E)\mid(a\Rightarrow p_1)
       \xrightarrow{1} (a,b\Rightarrow E).
\end{aligned}
\]
These exhaust the choices up to external permutation. The ordinary
root succedent stays $E$ throughout, so no splitting step occurs.
Consequently
\begin{equation}\label{def:bad}
 \mathcal D_{\mathrm{sp}}(U_E)=
 \{\,a,b\Rightarrow(a\cdot b)\wedge(b\cdot a)\,\}.
\end{equation}
Its only member has the semantic content $ab\le ab\wedge ba$.

Here is a finite countermodel.  Let $L=\{0,1,2,3\}$ with its usual order,
meet and join, with unit $e=3$, and with multiplication
\begin{equation}\label{def:mul}
 \begin{array}{c|rrrr}
 \cdot&0&1&2&3\\\hline
 0&0&0&0&0\\
 1&0&0&1&1\\
 2&0&0&2&2\\
 3&0&1&2&3
 \end{array}.
\end{equation}
The element $3$ is a two-sided identity, $0$ is absorbing, and every row
and column is nondecreasing.  Associativity can be checked structurally:
for $x,y\in\{0,1,2\}$, $xy=x$ if $y=2$, and $xy=0$ otherwise.  Thus for
$x,y,z$ in this subset, both associations equal $x$ if $y=z=2$, and
otherwise equal $0$.  Triples containing $3$ reduce by the unit law.

Define the residuals by
\[
 x\backslash z=\max\{y:xy\le z\},\qquad
 z/x=\max\{y:yx\le z\}.
\]
The sets are nonempty initial segments by monotonicity and absorption of
$0$, so these maxima exist and give both residuation equivalences.  For
completeness their tables are
\begin{equation}\label{def:residuals}
 \begin{array}{c|rrrr}
 x\backslash z&z=0&1&2&3\\\hline
 x=0&3&3&3&3\\
 1&1&3&3&3\\
 2&1&1&3&3\\
 3&0&1&2&3
 \end{array}
 \qquad
 \begin{array}{c|rrrr}
 z/x&z=0&1&2&3\\\hline
 x=0&3&3&3&3\\
 1&2&3&3&3\\
 2&0&1&3&3\\
 3&0&1&2&3
 \end{array}.
\end{equation}
Hence $L$ is a bounded integral residuated chain.  Since $e$ is greatest,
$xy\le e$ and $yx\le e$ always hold; in particular it satisfies the unit
weak-exchange condition.  Under $v(a)=1$, $v(b)=2$,
\[
 v(a\cdot b)=1,\qquad v(b\cdot a)=0,\qquad v(E)=0.
\]
The sequent in \eqref{def:bad} is therefore false.

Erase every numerical label in a putative proof, interpreting all
labelled occurrences by the same distinguished variable. The ordinary
rules, COM, WCM and WCT are sound on $L$, as are the generalized
lattice rules on a chain. After erasure the two copy blocks have the
same multiset of sequents. Deleting one of them therefore preserves
validity by ordinary external contraction. In Construction~3(ii) the
additional deleted background $G''_2$ is a copy of the retained $G'_2$,
so this deletion preserves validity for the same reason. This argument
also proves soundness when the two deletion rules are admitted without
their scheduling restrictions, for any set of numerical labels used in
the proof. The label-free sequent in~\eqref{def:bad} is false on $L$.
It is therefore unprovable even in this enlarged calculus, and hence
unprovable in the published calculus. The eligibility of the original
contraction~\eqref{def:contract} was verified directly above.
\end{proof}

The algebra \eqref{def:mul} is our countermodel for this auxiliary claim;
it is not the four-element algebra, with unit $2$, used for the different
cut-elimination obstruction in Lemma~1 of \cite{Wang2019}.

Definition~9 of~\cite{Wang2019} defines $\mathcal D$ on inputs carrying
a $\mathrm{GL}^{\Omega}_{\mathrm{cf}}$ proof. The contraction in
\eqref{def:contract} supplies a $\mathrm{GpsUL}^{\Omega}$ proof of $U_E$,
whereas Lemma~5 applies $\mathcal D$ to that conclusion. With the
domain of Definition~9 left unchanged, this inference lacks the required
domain justification. With the syntactic extension of its splicing
clause, the image is exactly~\eqref{def:bad}.

The uncontracted image explains the obstruction:
\begin{equation}\label{def:good}
 \mathcal D(H)=
 \{\,a,b\Rightarrow E\mid b,a\Rightarrow E\,\}.
\end{equation}
It is valid on every chain because one of $ab$ and $ba$ is their minimum.
Thus neither ordinary contraction after erasing labels nor
contraction-free splicing is at fault.  Rather, the two labelled blocks
have different ordered splicing images, and the stated compatibility
claim removes an essential alternative from \eqref{def:good}.

\subsection{Coherent copies and the proof obligations for density elimination}
\label{def:repair-scope}

For closed blocks $B,B'$ write $J,J'$ for their respective active label
sets. A coherent copy has a bijection $h:J\to J'$ and a matching $\pi$
of component occurrences that preserves ordinary formulas, succedent
positions and antecedent positions, and realizes $h$ on both sides.
For copies to be contracted we also require $J\cap J'=\varnothing$.
Such a bijection extends to a permutation of the ambient names by
interchanging each $j\in J$ with $h(j)\in J'$ and fixing names outside
the two supports.

Now fix a particular matching $\pi$ witnessing a published copy by
bijections $\sigma_l,\sigma_r$. A single bijection realizes this
\emph{same} matching if and only if
$\sigma_l(j)=\sigma_r(j)$ for every $j\in J$. Indeed, closedness gives
one antecedent and one succedent occurrence of each $j$. Their target
positions under $\pi$ force the values $\sigma_l(j)$ and $\sigma_r(j)$,
respectively. If the two values agree for each $j$, take $h=\sigma_l$.
The existence of a coherent copy allows a different matching of the
component occurrences. In~\eqref{def:U}--\eqref{def:V}, however, the
ordinary root is unique and the two producer antecedents $a,b$ are
distinct. Hence any matching must send the root to the root and match
the producers by their antecedents. The two root positions force
$1\mapsto3,\ 2\mapsto4$, while the producer succedents force
$1\mapsto4,\ 2\mapsto3$. Thus~\eqref{def:contract} is excluded by
coherent copying.

This change repairs the identification of producers with ordered consumer
positions, but does not by itself prove the original Lemma~7.  Two further
obligations are distinct.  First, the density-compatibility invariant must
be iterable through every auxiliary inference and every permitted
contraction.  The contraction-free provability assumption initializing a
single application of Lemma~5 cannot simply be reused after an earlier
contraction.  In splitting units one must also account for the choices of
left and right outputs and their backgrounds; matching labels alone does
not identify independently chosen proof-dependent partitions.  Second,
the global separation construction must generate the coherent blocks and
preserve closed boundaries, ordered slot positions and genuine contraction
histories.  In particular, the conclusion of the proof of Property~(A)
in the original Lemma~7 establishes contractibility of two copies; a
normalization step preserving the other invariants is needed to infer the
stated ``at most one copy'' property of the constructed object.

The replacement developed below addresses these as separate local and
global statements.  It uses the Cartesian family of density outputs for
the local invariant, and an ordered occurrence-based separation argument
for the global construction.  It does not claim that merely changing
Definition~5(iv) validates all contractions of the original algorithm.

Finally, Proposition~\ref{def:counter} concerns a universal intermediate
claim with all its stated premises satisfied.  It does not show that the
particular mismatched copies arise in every execution of the intended
separation algorithm, and it does not refute the conclusions of the
original density-elimination or standard-completeness theorems.  Their
published proof uses the defective compatibility claim, so a replacement
argument is required.  The original cut-elimination obstruction and the
contraction-free ordered splitting and splicing constructions are not
refuted by this counterexample.

\section{A Cartesian density transformation}\label{sec:cartesian-density}

We give a local proof transformation that is independent of any
separation algorithm. Its hypotheses and conclusion concern finite
families of ordinary hypersequents. In particular, it does not use a
density rule to establish its own admissibility.

\subsection{The auxiliary calculus}

Throughout this section, the target calculus is
$\mathsf K=\Gtwo$, with the rules in Section~\ref{sec:calculus}.
The auxiliary calculus below preserves the order of antecedent words.

Fix a distinguished variable $p$. In the auxiliary calculus, its
occurrences are replaced by labels $p_i$. A label occurs only as a whole
antecedent formula or as the whole succedent; it never occurs inside a
logical connective. An \emph{ordinary} formula contains neither $p$ nor any of these
labels. A labelled hypersequent is \emph{closed} when every label occurring
in it has exactly one antecedent occurrence and exactly one succedent
occurrence. Two such hypersequents are \emph{label-disjoint} when their
label sets are disjoint. Component and formula occurrences, rather than
just their displayed strings, are retained. The unlabelled skeleton is
obtained by replacing every $p_i$ by $p$, preserving all antecedent
positions and component multiplicities.

\begin{definition}[The linear labelled calculi]\label{def:cart-auxiliary}
The calculus $\mathsf A_0$ has the following rules and restrictions.
\begin{enumerate}
\item Its initial sequents are the ordinary initial sequents of
$\mathsf K$, with ordinary contexts, together with $p_i\Rightarrow p_i$.
\item It has communication, weak exchange, unit cut, unit-left, and the
logical rules of $\mathsf K$, with ordinary principal formulas. However,
ordinary meet-right and join-left are replaced, respectively, by
\begin{equation}\label{eq:cart-generalized-meet}
 \frac{G_1\mid\Gamma_1\Rightarrow A\quad
       G_2\mid\Gamma_2\Rightarrow B}
      {G_1\mid G_2\mid\Gamma_1\Rightarrow A\wedge B\mid
                         \Gamma_2\Rightarrow A\wedge B},
\end{equation}
\begin{equation}\label{eq:cart-generalized-join}
 \frac{G_1\mid\Gamma_1,A,\Lambda_1\Rightarrow C_1\quad
       G_2\mid\Gamma_2,B,\Lambda_2\Rightarrow C_2}
      {G_1\mid G_2\mid\Gamma_1,A\vee B,\Lambda_1\Rightarrow C_1\mid
                         \Gamma_2,A\vee B,\Lambda_2\Rightarrow C_2}.
\end{equation}
Here $A,B$ are ordinary; $C_1,C_2$ may be labelled succedents.
\item Every premise and conclusion is closed. The two premises of each
binary inference are label-disjoint. There is no external weakening,
external contraction, or unrestricted cut.
\end{enumerate}
Thus $\mathsf A_0$ is the contraction-free labelled calculus
$\mathrm{GL}^{\Omega}_{\mathrm{cf}}$, specified here without referring
to a proof-dependent density operation.

The calculus $\mathsf A_{\mathrm{ord}}$ extends $\mathsf A_0$ by
\begin{equation}\label{eq:cart-coherent-ec}
 \frac{G\mid B_1\mid B_2}{G\mid B_1}\ (\mathrm{EC}_{\mathrm c}),
\end{equation}
where $B_1,B_2$ are nonempty disjoint closed blocks and a single bijection
of labels takes $B_1$ to $B_2$, simultaneously on the two sides of every
sequent and preserving all ordinary formulas and ordered antecedent
positions. We call them \emph{coherent copies}. Zero-label blocks are
allowed. Identity steps on a whole hypersequent may also be used.
\end{definition}

The two generalized logical rules are derivable in $\mathsf K$.
For~\eqref{eq:cart-generalized-meet}, communication on the two premises
first gives the crossed components
$\Gamma_2\Rightarrow A\mid\Gamma_1\Rightarrow B$.
Ordinary meet-right with the first original premise forms the meet at
$\Gamma_1$; ordinary meet-right with the second forms it at $\Gamma_2$.
External contraction removes repeated backgrounds.
For~\eqref{eq:cart-generalized-join}, communication exchanges the
displayed occurrences of $A$ and $B$. The crossed components are
$\Gamma_1,B,\Lambda_1\Rightarrow C_1$ and
$\Gamma_2,A,\Lambda_2\Rightarrow C_2$.
Apply ordinary join-left first in the first context and then in the
second, and contract repeated backgrounds. These derivations use no
internal exchange.

We will also use \emph{unit insertion}, obtained by unit-left followed
by unit cut:
\begin{equation}\label{eq:cart-insertion}
 \frac{G\mid\Gamma,\Lambda\Rightarrow A\quad
       J\mid\Pi\Rightarrow t}
      {G\mid J\mid\Gamma,\Pi,\Lambda\Rightarrow A}.
\end{equation}
In particular, two words with succedent $t$ can be concatenated. All
uses of insertion below are finite derivations in $\mathsf K$.

\subsection{Closed units and ordered expansion}

For a closed labelled hypersequent $H$, form a directed graph whose
vertices are its component occurrences. The component with succedent
$p_i$ has one outgoing edge to the component containing the antecedent
occurrence of $p_i$. An ordinary-succedent component has no outgoing
edge. Closed subhypersequents are exactly unions of weakly connected
components of this graph: closure requires inclusion of both ends of
every edge, and the converse is immediate. Consequently its weakly
connected components are precisely the minimal nonempty closed units.

A finite weakly connected graph of this form has either one vertex of
outdegree zero and no cycle, or one directed cycle and no vertex of
outdegree zero. Indeed, following the unique outgoing edge reaches a
sink or a cycle; different such destinations cannot lie in one weak
component. Call these two kinds \emph{splicing units} and
\emph{splitting units}, respectively.

An \emph{actual splice} contracts the edge of a label $i$ whose producer
and consumer are distinct component occurrences. Up to external
permutation it replaces
\[
 U\Rightarrow p_i\ \mid\ A,p_i,B\Rightarrow D
 \quad\text{by}\quad A,U,B\Rightarrow D .
\]
All other occurrences are retained. In particular, a splice preserves
closedness, concatenates words in the displayed order and deletes exactly
one label. A \emph{complete splice sequence} ends when all remaining
label edges are loops at individual vertices. A sequence of length zero
is allowed. Each such sequence, together with its terminal hypersequent,
will be called an \emph{expansion witness}.

There is also a readout directly from the original graph. In a
splicing unit $U$, read the word at its ordinary root by recursive
ordered substitution: at each antecedent label read the word at its
unique producer in that position. The resulting ordinary sequent is
$s(U)$. In a splitting unit $C$, prescribe an original cycle label
$q_C$. Cut its edge, leaving the vertex producing $q_C$ as the root
and the occurrence consuming $q_C$ as a distinguished leaf.
The remaining edges form a rooted tree. The same ordered substitution,
retaining that leaf, gives the single component
\[
 L_{C,q_C},q_C,R_{C,q_C}\Rightarrow q_C .
\]
Here and below a cut of a graph edge retains both boundary occurrences
of its label. It is a choice in reading a cycle, independent of the cut
rule of the calculus.

\begin{definition}[Cartesian density family]\label{def:cartesian-family}
Let $\mathcal U(H)$ and $\mathcal C(H)$ denote the splicing and splitting
units of a closed $H$. For a list $\mathbf q=(q_C)_{C\in\mathcal C(H)}$
of prescribed cycle cuts, put
\[
 e_{C,q_C}(L)=(L_{C,q_C}\Rightarrow t),\qquad
 e_{C,q_C}(R)=(R_{C,q_C}\Rightarrow t).
\]
For each $\sigma:\mathcal C(H)\to\{L,R\}$ define
\[
 \Delta_{\mathbf q,\sigma}(H)
   =\mathop{\mid}_{U\in\mathcal U(H)}s(U)\ \mid\!
       \mathop{\mid}_{C\in\mathcal C(H)}e_{C,q_C}(\sigma(C)).
\]
Write $\Delta_{\mathbf q}(H)$ for the family indexed by these
$2^{|\mathcal C(H)|}$ endpoint choices. Occurrence multiplicities are
retained even when two displayed sequents coincide. An expansion witness
has its own family: take every ordinary terminal root and choose a left
or right endpoint at each terminal loop. Write $\mathsf{All}(H)$ when
every member of every expansion witness of $H$ has a $\mathsf K$ proof.
When there are no cyclic units there is just one endpoint choice.
\end{definition}

\begin{lemma}[Realization of cycle cuts and cycle-cut independence]
\label{lem:fss-cycle-cut}
For any prescribed $\mathbf q$ there is one expansion witness whose
member for each $\sigma$ is $\Delta_{\mathbf q,\sigma}(H)$ up to
external permutation. Conversely, every
expansion witness reads precisely this family for the original cycle
labels surviving in that witness, up to external permutation.
For any two prescribed $\mathbf q,\mathbf q'$ and the \emph{same}
choice $\sigma$ on the original units, the two displayed members are
interderivable in $\mathsf K$. Consequently, for each fixed $\mathbf q$,
\[
 \mathsf{All}(H)\quad\Longleftrightarrow\quad
   (\forall\sigma)\ \mathsf K\vdash\Delta_{\mathbf q,\sigma}(H).
\]
\end{lemma}
\begin{proof}
In each rooted or cut tree, contract its remaining tree edges, for
example starting at the leaves. None of these contractions is a loop:
contracting an edge of a tree gives another tree. Each contraction
substitutes exactly the word that the recursive readout assigns at
that occurrence. More explicitly, evaluating the remaining child words
after splicing $U$ into $A,p_i,B$ gives the concatenation of the
evaluated $A,U,B$; evaluating the children first gives the same
concatenation. A distinguished cut leaf is carried as one unchanged
letter in this identity. Thus the recursive readout is invariant at
each splice, including for empty substituted words. External permutations
only change the order in which components are listed. After the last
tree edge is contracted, each unit has its stated root or its
prescribed loop as its only component. This constructs the witness
before any endpoint is chosen.

For the converse, a splice contracts two vertices within one connected
unit and preserves that unit's connectedness. A terminal connected unit
therefore has exactly one vertex, since two distinct vertices in it
would require a nonloop edge. A unit with $n$ vertices had $n-1$ edges
in the rooted case and $n$ in the cyclic case. Exactly $n-1$ splices
occur within it, each deleting its own original edge. In the cyclic case
the surviving label belongs to the original cycle: every off-cycle edge
is a bridge of the underlying undirected multigraph, and contracting
other edges cannot identify the two ends of a bridge. Cut the surviving
edge in the original unit. Every splice in the given sequence now
contracts a tree edge and preserves the readout by the preceding
substitution identity. This identifies the two endpoint words, with
their order, for this particular witness.

To compare two prescribed cuts, first evaluate all branches attached to
the cycle. Its components can then be written
$\alpha_i,p_{i-1},\beta_i\Rightarrow p_i$ for $1\leq i\leq m$, with
indices taken modulo $m$. The cut at $p_k$ gives
\[
 L_k=\alpha_k\alpha_{k-1}\cdots\alpha_{k+1},\qquad
 R_k=\beta_{k+1}\beta_{k+2}\cdots\beta_k ,
\]
each word containing exactly $m$ indicated blocks. Changing $k$
cyclically rotates each of the two words separately. Whichever endpoint
$\sigma$ selects, one instance of weak exchange takes
$X,Y\Rightarrow t$ to $Y,X\Rightarrow t$ under the full external
background; the reverse instance restores the original word.
Apply this at the selected component of each affected unit. Finally,
the prescribed witness gives the forward implication in the display.
For the reverse implication take an arbitrary witness, identify its
surviving labels as above and transport the proof of the corresponding
fixed-cut member by those rotations.
\end{proof}

\begin{lemma}[Ordered one-hole locality]\label{lem:fss-locality}
Mark a consecutive, possibly empty, antecedent interval $P$ in a
component $S$ of a closed unit. In a cyclic unit cut at the succedent
label of $S$ if $S$ lies on the cycle, and at any cycle label otherwise.
In the readout of the resulting tree the expanded interval $\widehat P$
occupies one specified consecutive position. It contains the retained
placeholder if and only if $S$ lies on the cycle and $P$ contains the
antecedent occurrence supplied by the preceding cycle vertex.

For two label-disjoint input units of a communication step, exchange
their marked intervals. A readout of the resulting units is obtained
by exchanging the two expanded intervals in the root or cut-tree words,
then splicing the remaining edges between distinct root components.
Other units and their readouts are unchanged.
\end{lemma}
\begin{proof}
In a rooted tree, each antecedent position is either an ordinary letter,
the distinguished cut leaf, or the attachment of a unique child tree.
Write the word at $S$ as $A,P,B$ and detach the ordered list of letters
and child trees attached at the positions of $P$. Replace this list by
a marked gap. Following the unique path from $S$ to the root and
substituting all other child words gives a word $X,* ,Y$ with exactly
one marked gap: at each ancestor the unique child on this path is
inserted once in its ordered position. Reinstating the detached list
inserts its concatenated readout $\widehat P$ at that gap. For $P$
empty, the two boundaries select the original gap even when several
empty words have the same display.

In a rooted unit no cut leaf occurs. In a cyclic unit, if $S$ is off
the cycle, its descendant trees contain no cycle vertex
and hence no distinguished cut leaf. If $S$ is on the cycle, our cut
makes it the root. Exactly one antecedent position there leads to the
cut leaf: it is the position supplied by the preceding cycle vertex,
or the cut leaf itself for a one-vertex cycle. The detached list
contains the leaf precisely when it includes that position. This proves
both directions, also when the selected interval is empty.

For communication perform the detachment in both cut trees, using
distinct marked gaps. Each detached list contains whole child trees;
its active vertex stays in the complementary tree. Grafting the
opposite list into each gap therefore leaves two rooted ordered trees
with the original root succedents. At the level of the labelled
sequents, this exchanges precisely the two intervals: a transferred
label occurrence takes its entire supplying subtree with it, and its
producer and the consumer occurrences elsewhere in that subtree are
unchanged. Distinct label sets prevent any further identification.
For each original cut, its producer label still names its original root;
its marked leaf may have moved into the other tree. Contracting all
noncut edges now reads the two grafted trees. Reconnecting the retained
labels gives exactly these root components with the corresponding leaves
as consumer occurrences. Splice any such label between distinct
components until only ordinary roots and self-loops remain. This is
an actual expansion witness of the communicated hypersequent. Any
other choice of final cycle anchors gives corresponding members by
Lemma~\ref{lem:fss-cycle-cut}.
\end{proof}

For brevity, $\Delta(H)$ in the compatibility proofs below denotes the
family for any one prescribed list of cuts. Lemma~\ref{lem:fss-cycle-cut}
allows us to use cuts convenient for the rule and then to recover every
member of any prescribed target family, or of any actual expansion
witness. An empty interval in a loop word retains its marked gap on the
indicated side of the placeholder.

\subsection{Local compatibility}

\begin{lemma}[Coherent closed-block contraction]\label{lem:fss-cart-ec}
For an instance of~\eqref{eq:cart-coherent-ec}, proofs in $\mathsf K$
of all members of $\Delta(G\mid B_1\mid B_2)$ yield proofs of all
members of $\Delta(G\mid B_1)$.
\end{lemma}
\begin{proof}
A closed block is a union of minimal units. The common bijection matches
the splicing units of $B_1,B_2$ with identical ordered root expansions,
and matches their splitting units with identical endpoint words when
corresponding cuts are used. Such cuts may be assumed by
Lemma~\ref{lem:fss-cycle-cut}. Fix a choice for $G\mid B_1$ and choose
the same endpoint for each corresponding splitting unit in $B_2$.
In this diagonal member of the premise family, every component
contributed by $B_2$ duplicates the corresponding component from $B_1$.
Ordinary external contraction removes these duplicates and yields the
chosen conclusion member.
\end{proof}

In the remaining proofs, fix arbitrary choices in units untouched by
the rule and denote their contributions in the two premises by $E_1,E_2$.
These are whole external backgrounds. External weakening adds missing
passive components; external contraction removes repeated ones. Neither
operation changes an ordered antecedent.

\begin{lemma}[Communication compatibility]\label{lem:fss-cart-com}
If every Cartesian density member of each premise of a communication
instance in $\mathsf A_0$ has a $\mathsf K$ proof, then every Cartesian
density member of its conclusion has a $\mathsf K$ proof.
\end{lemma}
\begin{proof}
Only the unit containing each active component can change. By
Lemma~\ref{lem:fss-locality}, each affected unit can be written as an
ordinary root word or a loop word with a consecutive exchanged interval.
Lemma~\ref{lem:fss-cycle-cut} permits the indicated anchors. Each
exchanged interval either contains its loop placeholder or is ordinary.
We consider all possibilities.

\emph{No exchanged interval contains a loop placeholder.}
The two affected units remain distinct and retain their kinds. An
ordinary root changes by the displayed exchange. A loop changes only
in the endpoint word containing the interval. Fix a target endpoint
choice for each loop. If any chosen endpoint is unaffected, its premise
clause survives as a submultiset of the target clause; external weakening
supplies the rest. Otherwise both selected words are affected. Apply
ordinary communication directly to those words, whose succedents are
either their ordinary root formulas or $t$. Its active conclusions are
exactly the selected target words, with background $E_1\mid E_2$.
This also handles two ordinary roots.

\emph{Exactly one exchanged interval contains a placeholder, and the
other unit is an ordinary root.} Write the input loop and root as
\[
 \Gamma_1,A_1,p_1,B_1,\Lambda_1\Rightarrow p_1,
 \qquad \Gamma_2,Z,\Lambda_2\Rightarrow C,
\]
where the exchanged intervals are $A_1,p_1,B_1$ and $Z$.
After communication and expansion of $p_1$, they form one ordinary root
\[
 \Gamma_2,A_1,\Gamma_1,Z,\Lambda_1,B_1,\Lambda_2\Rightarrow C.
\]
The loop premise supplies both
$E_1\mid\Gamma_1,A_1\Rightarrow t$ and
$E_1\mid B_1,\Lambda_1\Rightarrow t$.
Weak exchange gives $A_1\Gamma_1\Rightarrow t$ and
$\Lambda_1B_1\Rightarrow t$, with background $E_1$ in each.
Insert these words immediately before and after $Z$ in the root,
using~\eqref{eq:cart-insertion} twice. Contract the duplicated $E_1$.
The result is the required target member.

\emph{Exactly one exchanged interval contains a placeholder, and the
other unit is a loop.} Use the same first loop. If the ordinary exchanged
interval $Z$ lies to the left of $p_2$, write the second loop as
\[
 \Gamma_2,Z,\Sigma_2,p_2,\Lambda_2\Rightarrow p_2.
\]
After communication and expansion of $p_1$, the single resulting loop
has endpoint words
\[
 L=\Gamma_2 A_1\Gamma_1 Z\Lambda_1 B_1\Sigma_2,
 \qquad R=\Lambda_2.
\]
For its left choice, rotate the two endpoint proofs from the first
premise as in the previous case, and insert their words on the two
sides of $Z$ in
$E_2\mid\Gamma_2,Z,\Sigma_2\Rightarrow t$.
For its right choice, weaken $E_2\mid\Lambda_2\Rightarrow t$ by $E_1$.
If $Z$ lies on the right, write the second loop as
$\Gamma_2,p_2,\Sigma_2,Z,\Lambda_2\Rightarrow p_2$.
Now the left endpoint is unchanged and the right endpoint becomes
$\Sigma_2 A_1\Gamma_1 Z\Lambda_1 B_1\Lambda_2$.
The same two insertions prove that choice; the unchanged choice follows
by weakening.

\emph{Both exchanged intervals contain their placeholders.}
The input loops and exchanged intervals are
\[
 \Gamma_i,A_i,p_i,B_i,\Lambda_i\Rightarrow p_i,
 \qquad \Pi_i=A_i,p_i,B_i\quad(i=1,2).
\]
The communicated components are
\[
 \Gamma_1,A_2,p_2,B_2,\Lambda_1\Rightarrow p_1,
 \qquad
 \Gamma_2,A_1,p_1,B_1,\Lambda_2\Rightarrow p_2.
\]
Expansion of $p_2$ gives one loop with endpoint words
\[
 L=\Gamma_1 A_2\Gamma_2 A_1,
 \qquad R=B_1\Lambda_2 B_2\Lambda_1.
\]
For $L$, rotate the left input choices to
$A_1\Gamma_1\Rightarrow t$ and
$A_2\Gamma_2\Rightarrow t$, concatenate them by unit insertion, and
rotate the result to $\Gamma_1 A_2\Gamma_2 A_1\Rightarrow t$.
The backgrounds combine to $E_1\mid E_2$.
For $R$, rotate the right input choices to
$\Lambda_1B_1\Rightarrow t$ and
$\Lambda_2B_2\Rightarrow t$, concatenate, and rotate to
$B_1\Lambda_2B_2\Lambda_1\Rightarrow t$.

These cases exhaust root and loop combinations and both possible
positions of every ordinary exchanged interval. In particular, an
internal off-cycle focus introduces no additional case: its interval
remains ordinary and consecutive by Lemma~\ref{lem:fss-locality}.
All proofs are finite compositions of the rules of $\mathsf K$.
\end{proof}

\begin{lemma}[Compatibility of the other rules]\label{lem:fss-cart-logical}
Every remaining rule of $\mathsf A_0$ preserves provability in
$\mathsf K$ of all Cartesian density members.
\end{lemma}
\begin{proof}
Fix a target member and make matching choices in all passive units.

For an ordinary initial sequent, its sole density member is that initial
sequent. The two choices for $p_i\Rightarrow p_i$ are both
$\Rightarrow t$. The restricted bottom and top axioms contain no labels.

For a unary antecedent rule, its ordinary principal occurrence occupies
a unique interval in one root word or one endpoint word by
Lemma~\ref{lem:fss-locality}. If that word is selected, apply the same
rule in its expanded context. If the other endpoint of a loop is
selected, the selected sequent is unchanged. This handles product-left,
either meet-left rule, and unit-left. For a unary right rule, the active
ordinary succedent is the root of a splicing unit. Join-right is
immediate. In either residual-right rule, the discharged ordinary
formula remains at the correct outer end of the flattened word, so the
same right rule applies. Weak exchange has an ordinary root succedent
$t$: its displayed blocks expand to $X,Y$, and the target step is exactly
$X,Y\Rightarrow t$ to $Y,X\Rightarrow t$.

For product-right, the two active units are splicing units with root
words $X\Rightarrow A$ and $Y\Rightarrow B$. The conclusion merges
them into the root $X,Y\Rightarrow A\cdot B$, obtained by ordinary
product-right. Passive choices are unchanged.

For residual-left, the argument-producing premise has an ordinary root
$X\Rightarrow A$. If the selected word of the receiving unit contains
the principal occurrence $B$, write it as $L,B,R\Rightarrow D$.
The ordered residual-left steps on the flattened words are
\[
 \frac{E_1\mid L,B,R\Rightarrow D\quad E_2\mid X\Rightarrow A}
      {E_1\mid E_2\mid L,X,A\to B,R\Rightarrow D},
\]
\[
 \frac{E_2\mid X\Rightarrow A\quad E_1\mid L,B,R\Rightarrow D}
      {E_1\mid E_2\mid L,A\rightsquigarrow B,X,R\Rightarrow D}.
\]
Here $D$ is the ordinary root formula or the endpoint succedent $t$.
Locality identifies these words with the target choices. If the receiving
unit is a loop and the other endpoint is selected, that endpoint has
not changed. Its receiving-premise clause alone therefore yields the
target by external weakening. The argument-producing root has been
absorbed into the unselected endpoint and is not needed in that clause.

Unit cut is the same argument with the marked ordinary occurrence $t$
and argument root $X\Rightarrow t$. When the selected word contains
that occurrence, use
\[
 \frac{E_1\mid L,t,R\Rightarrow D\quad E_2\mid X\Rightarrow t}
      {E_1\mid E_2\mid L,X,R\Rightarrow D}.
\]
When it selects the other endpoint, use the unchanged receiving-premise
clause and external weakening.

For generalized meet-right, both active units have ordinary root
succedents. Apply the derivation of~\eqref{eq:cart-generalized-meet}
given above to their root words. For generalized join-left, if a target
choice avoids one principal occurrence, that entire premise clause
survives in the target, which follows by external weakening. Otherwise
the two selected words have forms
$L_1,A,R_1\Rightarrow D_1$ and $L_2,B,R_2\Rightarrow D_2$.
Apply the derivation of~\eqref{eq:cart-generalized-join} in these
expanded contexts. It gives exactly the two selected target words.

All binary steps may duplicate ordinary external backgrounds; external
contraction removes those copies. There are no other rules in
$\mathsf A_0$, apart from harmless identity steps on a hypersequent.
\end{proof}

\begin{theorem}[Strong local Cartesian density transformation]
\label{thm:fss-cartesian-local}\label{thm:cartesian-local}
Let a finite $\mathsf A_{\mathrm{ord}}$ derivation end in $H$ and have
closed hypothesis blocks $H_1,\ldots,H_m$. A leaf may use a complete
copy of the corresponding input with its labels uniformly freshened.
If $\mathsf{All}(H_i)$ holds for every input, then $\mathsf{All}(H)$.
Equivalently, prescribe one list of cycle cuts at each input and one
at the conclusion. Proofs of all endpoint choices for the prescribed
input families yield proofs of all endpoint choices for the prescribed
conclusion family. In particular, a closed derivation without
hypotheses yields proofs of the members of every actual expansion of $H$.
\end{theorem}
\begin{proof}
Induct on the derivation, with $\mathsf{All}$ as invariant. At a
hypothesis, a uniform injective renaming on its finite label set
is a bijection onto the labels of the copy. Transport a splice path
back along this bijection on occurrences and carry its endpoint choices
with it. The ordinary output is unchanged up to external permutation.
Initial sequents are covered by Lemma~\ref{lem:fss-cart-logical}.
At an inference take an arbitrary actual target expansion and a
choice in it. Lemma~\ref{lem:fss-cycle-cut} identifies its original
units and surviving cuts. Lemmas~\ref{lem:fss-cart-com},
\ref{lem:fss-cart-logical} and~\ref{lem:fss-cart-ec} supply a proof
for the corresponding target clause using the required premise members;
these are available by the induction hypotheses. Restoring the given
cuts uses weak exchange and external permutation. When a rule uses both
endpoints, apply the premise hypothesis to each choice to obtain the
two required proofs. This proves the
invariant at the conclusion. For the fixed-cut formulation, first
use Lemma~\ref{lem:fss-cycle-cut} at each open input to obtain
$\mathsf{All}(H_i)$; at the conclusion use the witness realizing the
prescribed cuts. Every step in this argument takes place in $\mathsf K$
or manipulates the finite labelled occurrence data.
\end{proof}

\section{Strong density elimination}\label{sec:strong-density}

Here \emph{strong} refers to the form of the density transformation with
finitely many producer and consumer components: each producer is paired with
each consumer in the conclusion. It does not refer to strong standard
completeness, which concerns consequence from arbitrary sets of assumptions.
The distinct arbitrary-theory completeness result is established in
Section~\ref{sec:strong-standard}.

\begin{theorem}[Strong density elimination]\label{thm:strong-density}
Let
\[
G_0=G\mid\{\Pi_i\Rightarrow p\}_{i=1}^{m}
\mid\{\Gamma_j,p,\Delta_j\Rightarrow A_j\}_{j=1}^{n},
\qquad m,n\ge1,
\]
where $p$ is absent from $G$, from each formula of $\Pi_i,\Gamma_j,\Delta_j$,
and from each $A_j$. The displayed families are indexed occurrences and
may contain repetitions. If $\G\vdash G_0$, then
\[
\G\vdash D_0(G_0):=
G\mid\{\Gamma_j,\Pi_i,\Delta_j\Rightarrow A_j:
          1\le i\le m,\ 1\le j\le n\}.
\]
\end{theorem}
\begin{proof}
By Theorem~\ref{uc:normalization}, choose a proof of $G_0$ whose
cuts are all on $t$. Each component of the original $G_0$ has at most
one occurrence of $p$: a background component has none, and every
producer or consumer has exactly its displayed occurrence. This is the
slot bound required before pruning in Proposition~\ref{op:preparation}
and Theorem~\ref{olg:source-free}. Apply those results to the normalized proof. This gives a nonempty closed theorem $H$
of $\Om$ consisting only of permitted decorations of terminal skeletons.
Forgetting the source, protection and copy-family annotations gives an
actual derivation in Definition~\ref{def:cart-auxiliary}: the annotated
contractions impose additional matching conditions and are instances of
the coherent closed-block contraction already justified there.
Apply the closed-proof case of Theorem~\ref{thm:cartesian-local}.
Every member of each actual expansion has a $\Gtwo$ proof and hence,
by inclusion of the unit-cut rules, a $\G$ proof.

Every component of $H$ has at most one live label occurrence. A cyclic
unit would require at least two incidences at each cycle vertex, counting
the two ends of a loop separately; hence no cyclic unit occurs. Every live
unit is a producer-consumer pair, and every other unit is an isolated
ordinary sequent. The Cartesian image is therefore a singleton, whose
unique hypersequent $K$ is provable in $\G$.

A live pair contributes exactly
$\Gamma_j,\Pi_i,\Delta_j\Rightarrow A_j$ for some $i,j$.
An unchanged $p$-free terminal component contributes a component of $G$.
A constant-filled consumer is
$\Gamma_j,\top,\Delta_j\Rightarrow A_j$.
Cut the initial sequent $\Pi_i\Rightarrow\top$ into this component to
obtain $\Gamma_j,\Pi_i,\Delta_j\Rightarrow A_j$, choosing any
$1\le i\le m$. A constant-filled producer is $\Pi_i\Rightarrow\bot$;
cut it into the initial sequent
$\Gamma_j,\bot,\Delta_j\Rightarrow A_j$, choosing any
$1\le j\le n$. These are ordinary cuts on $p$-free formulas in the target
calculus $\G$ and may be performed under the current external background.
The assumptions $m,n\ge1$ ensure that the choices exist even when some
terminal types disappeared during preprocessing.

After these replacements every component is a component of $D_0(G_0)$.
Use ordinary external contraction to remove excess multiplicities and
external weakening to add missing components, obtaining exactly
$D_0(G_0)$.
\end{proof}

\begin{corollary}[Elimination of the density rule]\label{cor:density-elimination}
Every $\mathrm{GpsUL}^{*}_{D}$ proof can be transformed into a
$\mathrm{GpsUL}^{*}$ proof with the same conclusion.
\end{corollary}
\begin{proof}
In a finite proof containing density choose an occurrence with no density
in its premise subproof. Its single characteristic variable is fresh for
the external background and all displayed contexts. The case $m=n=1$ of
Theorem~\ref{thm:strong-density} supplies a proof of its conclusion without
density. Replace the chosen occurrence and its premise subproof by this
proof. The number of density occurrences decreases, and iteration ends
with a density-free proof.
\end{proof}

\section{From density admissibility to standard completeness}
\label{sec:fss-semantics}

The dense countermodel and its completion give weak standard completeness
for the density-augmented calculus. The final corollary for the original
calculus uses Corollary~\ref{cor:density-elimination}. Write $\G_D=\G+(D)$, where $\G$ includes its
ordinary cut rule and
\[
 \frac{H\mid\Pi\Rightarrow p\mid\Gamma,p,\Lambda\Rightarrow B}
      {H\mid\Gamma,\Pi,\Lambda\Rightarrow B}\ (D),
 \qquad p\notin\operatorname{Var}(H,\Gamma,\Pi,\Lambda,B).
\]
A bounded residuated chain is a linearly ordered bounded lattice equipped
with an associative product, a unit $e$, and two residuals satisfying
\[
 xy\le z\quad\Longleftrightarrow\quad y\le x\backslash z
       \quad\Longleftrightarrow\quad x\le z/y.
\]
We interpret $t$ as $e$ and the lattice constants as the bounds. The connectives
$A\to B$ and $A\rightsquigarrow B$ denote respectively
$a\backslash b$ and $b/a$. The value of an empty antecedent is $e$.
A sequent $\Gamma\Rightarrow A$ holds when the ordered product of
$\Gamma$ is at most the value of $A$; a hypersequent holds when at least
one component holds. A formula is valid when its value is at least $e$.
A \emph{standard algebra} has lattice reduct $[0,1]$. Its product need
not be commutative, and no additional requirement that $e<1$ is imposed.
All completeness assertions in this section concern theorems, without premises.

\begin{lemma}\label{lem:fss-dense-weak-exchange}
Every dense bounded residuated chain satisfies
$ab\le e\Longrightarrow ba\le e$. Consequently every standard algebra
validates weak exchange.
\end{lemma}
\begin{proof}
We use the direct residuation argument of~\cite[Lemma~10.4]{Wang2026}.
Suppose $ab\le e<ba$. Then $b\backslash e<a$. By density choose
$b\backslash e<r<a$, so $br>e$. Associativity and monotonicity give
\[
             a=ae\le a(br)=(ab)r\le er=r<a,
\]
a contradiction. Interchanging $a,b$ gives unit cyclicity.
\end{proof}

\begin{lemma}[Countable dense countermodel]
\label{lem:fss-dense-countermodel}
If $\G_D\nvdash\Rightarrow F$, there is a countable dense bounded
residuated chain satisfying weak exchange, and a valuation in it with
$v(F)<e$.
\end{lemma}
\begin{proof}
Work in a countable language with infinitely many variables, and enumerate
all sequents as $S_0,S_1,\ldots$. For a finite set $B$ of sequents, its
associated hypersequent contains each member once. We construct increasing
finite sets $B_n$ whose associated hypersequents are underivable in
$\G_D$, starting with $B_0=\{\Rightarrow F\}$.

At stage $n$, adjoin $S_n$ if underivability is preserved, and otherwise
do not adjoin it; call the result $C_n$. If $S_n=A\Rightarrow B$ belongs
to $C_n$, choose a variable $p$ absent from $C_n$ and put
\begin{equation}\label{eq:fss-density-witness}
 B_{n+1}=C_n\cup\{A\Rightarrow p,\ p\Rightarrow B\}.
\end{equation}
This set is still underivable. Indeed, a proof of its hypersequent,
followed by $(D)$ and external contraction, would prove $C_n$, since
$A\Rightarrow B$ already occurs in $C_n$. If the condition fails, put
$B_{n+1}=C_n$.

Let $B_\infty=\bigcup_nB_n$. No finite nonempty hypersequent consisting
of members of $B_\infty$ is provable: a proof would yield a proof of some
$B_n$ by external contraction and weakening. For every
$S\notin B_\infty$ there is a finite hypersequent $H_S$, consisting of
members of $B_\infty$, such that
\begin{equation}\label{eq:fss-positive-certificate}
 \G_D\vdash H_S\mid S.
\end{equation}
It is the certificate for refusing $S$ at its enumeration stage. A refused
sequent cannot enter later, since its certificate would then be a
provable all-$B_\infty$ hypersequent. In particular, every eventual member
$A\Rightarrow B$ of $B_\infty$ receives the two witnesses at its own
enumeration stage.

Call a sequent \emph{positive} if it is outside $B_\infty$. Positive
sequents contain the initial sequents and are closed under the ordinary
single-conclusion bounded full Lambek rules, cut, and weak exchange.
To verify this, take certificates~\eqref{eq:fss-positive-certificate}
for the premises, apply the desired rule with their finite bad contexts,
and contract duplicated context components. A bad conclusion would give
a provable all-bad hypersequent. If a common context is needed, add it
by external weakening at the roots of the completed certificate proofs.
Nothing is inserted inside a density inference, so no eigenvariable
condition is disturbed.

Define $A\preccurlyeq B$ when $A\Rightarrow B$ is positive. Identity and
cut make this a preorder. Moreover, COM applied to two identity sequents
proves
\[
 A\Rightarrow B\mid B\Rightarrow A.
\]
The two components cannot both be bad, so $\preccurlyeq$ is total.
Here are the details needed for its algebraic quotient. The product rules
and cut identify positivity of $A,B\Rightarrow C$ with
$A\cdot B\preccurlyeq C$; multiplication is monotone by product-right.
The residual rules and cut give
\[
 A\cdot B\preccurlyeq C
 \ \Longleftrightarrow\ B\preccurlyeq A\to C
 \ \Longleftrightarrow\ A\preccurlyeq B\rightsquigarrow C.
\]
For example, the reverse direction of the first residual equivalence
uses cut with $A,A\to C\Rightarrow C$. These adjunctions imply the
required antitonicity and monotonicity of each residual. The lattice
rules give the universal properties of meet and join. Associativity
and the two unit laws follow by the product and unit rules; the bound
axioms give the least and greatest elements. Thus all operations respect
the symmetric part of $\preccurlyeq$, and the quotient $L$ of the formula
algebra is a countable bounded residuated chain. Closure of positive
sequents under WCM gives $xy\le e\Longrightarrow yx\le e$ in $L$.

Unit-left and cut with $\Rightarrow t$ show that $\Rightarrow A$ is
positive exactly when $t\Rightarrow A$ is positive. Since
$\Rightarrow F\in B_\infty$, the canonical valuation therefore satisfies
$[F]<[t]=e$. Finally, if $[A]>[B]$, then $A\Rightarrow B$ is bad. Its
witnesses in the construction are both bad, and totality gives
\[
 [A]>[p]>[B].
\]
The chain is consequently dense and refutes $F$.
\end{proof}

\begin{lemma}[Completion of a countable dense chain]
\label{lem:fss-standard-completion}
Every nontrivial countable dense bounded residuated chain embeds into a
standard algebra. The embedding preserves both residuals and all
designated constants.
\end{lemma}
\begin{proof}
By a back-and-forth enumeration, a countable dense bounded chain is
order-isomorphic to $\mathbb Q\cap[0,1]$: fix the endpoints first and
alternately insert the next unused element of either enumeration into
the required open interval. Identify the given chain with this dense
subset $L$ of $[0,1]$, transporting its operations. This order embedding
preserves every existing supremum. Indeed, a real upper bound strictly
below an embedded supremum would, by density, give a smaller upper bound
already in $L$.

For $x,y\in[0,1]$ define
\begin{equation}\label{eq:fss-lower-product}
 x*y=\sup\{ab:a,b\in L,\ a\le x,\ b\le y\}.
\end{equation}
Monotonicity shows that this extends the original product. We check
preservation of arbitrary joins in its first coordinate. Let $x=\sup X$.
The inequality $\sup_{u\in X}(u*y)\le x*y$ follows from monotonicity.
Conversely, fix $a,b\in L$ with $a\le x,b\le y$. If $a<x$, some $u\in X$
satisfies $a<u$, whence $ab\le u*y$. If $a=x>0$, density and residuation
in $L$ give the existing suprema
\[
 a=\sup\{c\in L:c<a\},\qquad
 ab=\sup\{cb:c\in L,c<a\}.
\]
They are preserved by the order embedding. Each $c<a=x$ is below some
$u\in X$, so $ab\le\sup_{u\in X}(u*y)$ again. The case $a=0$, including
the empty join case, uses absorption by the bottom. Taking the supremum
over $a,b$ proves the required equality. The proof in the second
coordinate is identical and does not exchange the order of factors.

Every $x\in[0,1]$ is the supremum of $L\cap[0,x]$. Separate join
preservation therefore expands $(x*y)*z$ and $x*(y*z)$ to the same
supremum of $(ab)c=a(bc)$ over $a\le x,b\le y,c\le z$ in $L$.
The product is associative, and its unit is the embedded $e$, since
$e*x=x*e=\sup(L\cap[0,x])=x$. Join preservation supplies the two residuals
as the right adjoints of the multiplication maps.

These residuals extend the original ones. For example, let
$r=a\backslash b$ in $L$. Then $a*r=ar\le b$. If $x>r$, choose
$c\in L$ with $r<c\le x$. Residuation in $L$ yields $ac>b$, hence
$a*x>b$. Thus the new left residual of $a,b$ is exactly $r$; the right
residual is treated in the same way. The lattice operations, bounds,
and designated constants are also preserved. This gives the asserted
embedding. The completed chain satisfies weak exchange by
Lemma~\ref{lem:fss-dense-weak-exchange}.
\end{proof}

\begin{theorem}[Density-augmented standard completeness]
\label{thm:fss-density-standard}
For every formula $F$, the following are equivalent:
\begin{enumerate}
\item $\G_D\vdash\Rightarrow F$;
\item $F$ is valid in every dense bounded residuated chain satisfying
weak exchange;
\item $F$ is valid in every standard algebra.
\end{enumerate}
\end{theorem}
\begin{proof}
All rules of $\G$ are sound in the specified chains. In particular,
COM is sound because the exchanged block values are comparable, and
the smaller one can replace the larger in one of the two active
premises. Weak exchange is sound by hypothesis, and in standard
algebras by Lemma~\ref{lem:fss-dense-weak-exchange}.

For $(D)$, suppose its conclusion is false under a valuation. Every
component of $H$ is false. Write $\pi,\gamma,\lambda,b$ for the values
of $\Pi,\Gamma,\Lambda,B$. Since $\gamma\pi\lambda>b$, residuation gives
$\pi>r:=\gamma\backslash(b/\lambda)$. Choose $r<c<\pi$ and assign the
fresh variable $p$ the value $c$. The components $\Pi\Rightarrow p$ and
$\Gamma,p,\Lambda\Rightarrow B$ are then both false, while the other
components are unchanged. A universally valid premise therefore has a
universally valid conclusion. This proves (1)$\Rightarrow$(2),(3).

If (1) fails, Lemma~\ref{lem:fss-dense-countermodel} supplies a dense
chain refuting $F$, disproving (2). That chain is nontrivial because it
refutes a formula. Lemma~\ref{lem:fss-standard-completion} embeds it
into a standard algebra, preserving the valuation of every formula;
thus (3) fails as well.
\end{proof}

\begin{corollary}[Weak standard completeness]
\label{cor:fss-weak-standard}
For every formula $F$,
\[
 \G\vdash\Rightarrow F
 \quad\Longleftrightarrow\quad
 F\text{ is valid in every standard algebra}.
\]
Thus the formula-theorem standard-completeness conclusion is retained.
\end{corollary}
\begin{proof}
Standard validity gives a $\G_D$ proof by
Theorem~\ref{thm:fss-density-standard}. Corollary~\ref{cor:density-elimination}
transforms that finite derivation into a $\G$ proof. The converse is soundness.
\end{proof}

This density argument uses no density principle relative to arbitrary
theories and only completes the dense countermodels constructed above.
The next section supplies a separate arbitrary-theory argument, using
relative chain countermodels and a full bounded standard representation
that also applies to nondense countable chains. The Hilbert formulation
uses the theorem correspondence proved in Appendix~\ref{sec:hilbert-explicit}.

\section{Strong standard completeness for arbitrary theories}
\label{sec:strong-standard}

We now consider consequence from arbitrary sets of formulas in the
independent-parameter calculus $\mathrm H_{\mathrm{sl}}^{*}$ of
Appendix~\ref{sec:hilbert-explicit}. A corresponding arbitrary-theory
standard completeness statement is already given
in~\cite[Corollary~10.7]{Wang2026}. Here we prove the required relative
chain-countermodel result directly for the explicit calculus and use the
external bounded embedding for the final transfer. The passage from
theoremhood to this consequence relation requires a separate argument. We construct a chain
countermodel satisfying an entire theory under one valuation, and then
transport that valuation through a full bounded standard embedding.

Fix a countable set $V$ of propositional variables. The language contains
the constants $t,\bot,\top$ and may also contain the independent constant
$f$. We reserve $t$ for the formula denoting the monoidal unit and $e$
for the corresponding algebraic element. Thus $v(t)=e$, whereas $v(f)$
is an arbitrary distinguished element, independent of $e,\bot,\top$.
Substitution for ordinary variables fixes $f$. For any $V$, encode the
atoms of the pointed language by $\operatorname{Option}(V)$, using
$\operatorname{some}(p)$ for $p\in V$ and $\operatorname{none}$ for $f$.
If the underlying standard algebra is $S$, a valuation of these atoms
is equivalently a pair $(v,c)$ with $v:V\to S$ and $c\in S$: it
takes $\operatorname{some}(p)$ to $v(p)$ and $\operatorname{none}$ to $c$.
Fixing the pointed algebra $(S,c)$ fixes $c$ while $v$ varies.
The theory considered below may be an arbitrary set of formulas.

For a theory $\Gamma$, the notation
$\Gamma\vdash_{\mathrm H_{\mathrm{sl}}^{*}}A$ means that there is a finite
derivation using the displayed axiom schemata, the four Hilbert rules,
and specified members of $\Gamma$ as additional leaves. Such a leaf is
its stated formula, not every substitution instance of that formula.
Repeated uses of a premise are allowed. We abbreviate this relation by
$\Gamma\vdash A$ within this section.

Let $\mathsf{Std}$ denote the class of bounded residuated algebras on
the usual real unit interval, with their actual multiplication, both
residuals, and monoidal unit. The unit need not be an endpoint; the case
$e=1$ is allowed. These algebras satisfy weak exchange by
Lemma~\ref{lem:fss-dense-weak-exchange}. When $f$ is present, the class
includes all choices of its interpretation. Write
\[
 v\models_S\Gamma\quad\Longleftrightarrow\quad
 e_S\le v(A)\quad\text{for every }A\in\Gamma,
\]
and define
\[
 \Gamma\models_{\mathsf{Std}}F
 \quad\Longleftrightarrow\quad
 \text{for every $S\in\mathsf{Std}$ and every valuation $v$, }
 v\models_S\Gamma\ \Longrightarrow\ e_S\le v(F).
\]
In the pointed language the two outer quantifiers range over all
pointed algebras $(S,c)$ and all ordinary valuations $v:V\to S$.
The same $S,c,v$ interpret every member of $\Gamma$ and the conclusion
$F$. The Hilbert calculus for this consequence relation has precisely
the four rules of Appendix~\ref{sec:hilbert-explicit}.

\subsection{Finite support and the algebra relative to a theory}

\begin{lemma}[Finite support of actual derivations]
\label{lem:theory-finite-support}
For every theory $\Gamma$ and formula $F$,
\[
 \Gamma\vdash F
 \quad\Longleftrightarrow\quad
 \text{there are $A_1,\ldots,A_n\in\Gamma$ such that }
 \{A_1,\ldots,A_n\}\vdash F,
\]
where $n$ may be zero. More precisely, every supplied derivation can be
rebuilt over the finite list of its actual assumption leaves.
\end{lemma}
\begin{proof}
Traverse the given finite tree. An assumption leaf contributes its
formula, an axiom leaf contributes the empty list, a modus ponens node
concatenates the lists of its two premises, and each unary node keeps
the list of its premise. Every entry belongs to $\Gamma$. Induction on
the same tree reconstructs a derivation over the membership set of this
list: at modus ponens, weaken each reconstructed premise to the
membership set of the concatenated list and replay the rule. The other
constructors are replayed directly. This construction retains repeated
uses as repeated list entries. Conversely, a proof over a list of
members of $\Gamma$ is a proof over $\Gamma$ by weakening.
\end{proof}

For a fixed $\Gamma$, put
\[
 A\preceq_\Gamma B\quad\Longleftrightarrow\quad
 \Gamma\vdash A\to B,
 \qquad
 A\equiv_\Gamma B\quad\Longleftrightarrow\quad
 A\preceq_\Gamma B\text{ and }B\preceq_\Gamma A.
\]
Denote the quotient by $L_\Gamma$, and the class of $A$ by $[A]_\Gamma$.

\begin{lemma}[Relative Lindenbaum algebra]
\label{lem:relative-lindenbaum}
The relation $\equiv_\Gamma$ is a congruence, and the operations induced
by the formula constructors make $L_\Gamma$ a bounded residuated
lattice. Its unit is $e_\Gamma=[t]_\Gamma$, its order is
\[
 [A]_\Gamma\le[B]_\Gamma
 \quad\Longleftrightarrow\quad\Gamma\vdash A\to B,
\]
and it satisfies
\begin{equation}\label{eq:strong-designated}
 e_\Gamma\le[A]_\Gamma\quad\Longleftrightarrow\quad\Gamma\vdash A.
\end{equation}
It is unit-cyclic: $ab\le e_\Gamma$ if and only if $ba\le e_\Gamma$.
The canonical valuation sends every formula to its own class, including
$f$ when present. No totality of this quotient is asserted at this stage.
\end{lemma}
\begin{proof}
All derivations in this proof are over the same $\Gamma$. A1 gives
reflexivity of $\preceq_\Gamma$, and A2 followed by two applications of
modus ponens gives transitivity. If $\Gamma\vdash A$ and
$\Gamma\vdash B$, then A5 gives $\Gamma\vdash A\cdot B$.
Applying $(\mathrm{ADJ}_{\mathrm U})$ to both premises and then A8 gives
$\Gamma\vdash A\wedge B$. Thus two already available consequences can
be conjoined without changing the theory.

The derived currying and uncurrying constructions give
\begin{align*}
 \Gamma\vdash A\cdot B\to C
 &\quad\Longleftrightarrow\quad
 \Gamma\vdash B\to(A\to C),\\
 \Gamma\vdash A\cdot B\to C
 &\quad\Longleftrightarrow\quad
 \Gamma\vdash A\to(B\rightsquigarrow C).
\end{align*}
For the first equivalence, compose A5 with the instance of A2 induced
by $A\cdot B\preceq_\Gamma C$. The converse follows by modus ponens with A6. For the second, start with $B\to(A\to C)$, convert its outer
arrow to $\rightsquigarrow$ by right normalization and A7, and apply A4.
Conversely, uncurrying A3 gives
$(B\rightsquigarrow C)\cdot B\preceq_\Gamma C$. If
$A\preceq_\Gamma B\rightsquigarrow C$, right multiplication by $B$
and transitivity give $A\cdot B\preceq_\Gamma C$. Here monotonicity
of right multiplication already follows from the first equivalence.
If $X\preceq_\Gamma Y$, then
$B\preceq_\Gamma Y\to(Y\cdot B)\preceq_\Gamma X\to(Y\cdot B)$
by A5 and antecedent antitonicity, so $X\cdot B\preceq_\Gamma Y\cdot B$.
The antitonicity step uses A2 with its inner antecedent discharged
by the proof of $X\to Y$. In general, given proofs of $P\to(Q\to R)$
and $Q$, compose $P\to Q\cdot P$ (left normalization) with
$Q\cdot P\to R$ (A6). All these proofs use the same $\Gamma$.

The two equivalences give the remaining monotonicity of product and the
variance of both residuals. A9--A14 and the conjunction construction give
lattice monotonicity. Thus all operations respect $\equiv_\Gamma$ and the
two adjunctions descend to the quotient.

For associativity, first observe that in the quotient
$(xy)\backslash w=y\backslash(x\backslash w)$. The inequality from
right to left is A6. For the other direction A2 gives
\[
 (xy)\backslash w
 \le (x\backslash(xy))\backslash(x\backslash w).
\]
Since $y\le x\backslash(xy)$ by A5, antecedent antitonicity
bounds the right side by $y\backslash(x\backslash w)$.
Successive adjunctions now give $(xy)z\le w$ if and only if
$x(yz)\le w$ for every $w$, hence $(xy)z=x(yz)$.
For the unit laws, uncurrying A16 gives $t\cdot A\preceq_\Gamma A$.
Applying left normalization to A15 gives the converse.
A15 and A5 give $A\preceq_\Gamma A\cdot t$. Applying A16 to the
proved formula $A\to A$ and uncurrying gives $A\cdot t\preceq_\Gamma A$.
Finally A9--A14 give the universal properties of meet and join, and
A17 and A18 give the bounds.

For~\eqref{eq:strong-designated}, its left side means
$\Gamma\vdash t\to A$; modus ponens with A15 yields $A$.
Conversely, apply modus ponens to a proof of $A$ and A16,
$A\to(t\to A)$. In particular, this equivalence uses no deduction
theorem for an arbitrary set of premises.

WCM yields $e_\Gamma/a\le a\backslash e_\Gamma$. If
$ab\le e_\Gamma$, residuation gives
$a\le e_\Gamma/b\le b\backslash e_\Gamma$, hence $ba\le e_\Gamma$.
Interchanging $a,b$ gives the converse. Finally, evaluation at
$p\mapsto[p]_\Gamma$, with $f$ interpreted by $[f]_\Gamma$, sends
each formula to its class by induction on its construction.
\end{proof}

For reference, Table~\ref{tab:relative-derived} collects the derived
constructions just used. Each statement concerns derivations over the same
fixed theory $\Gamma$; none invokes a deduction theorem discharging an
arbitrary premise. In the third row, a proof of $Q$ yields
$P\to Q\cdot P$ by $(\mathrm{PN}_{\to})$, while A6 converts
$P\to(Q\to R)$ into $Q\cdot P\to R$. Transitivity then gives the stated
conclusion, retaining the order of the factors.

\begin{table}[htbp]
\centering\small
\caption{Derived constructions for the relative Lindenbaum algebra.}
\label{tab:relative-derived}
\begin{tabularx}{\textwidth}{@{}p{0.21\textwidth}X@{}}
\toprule
Construction & Derived statement over $\Gamma$\\
\midrule
Left currying & $\Gamma\vdash A\cdot B\to C$ iff
$\Gamma\vdash B\to(A\to C)$\\[3pt]
Right currying & $\Gamma\vdash A\cdot B\to C$ iff
$\Gamma\vdash A\to(B\rightsquigarrow C)$\\[3pt]
Insertion of a proved formula & From $\Gamma\vdash P\to(Q\to R)$ and
$\Gamma\vdash Q$, derive $\Gamma\vdash P\to R$\\[3pt]
Product monotonicity & $X\preceq_\Gamma Y$ implies
$X\cdot B\preceq_\Gamma Y\cdot B$ and
$B\cdot X\preceq_\Gamma B\cdot Y$\\[3pt]
Designation & $\Gamma\vdash A$ iff $\Gamma\vdash t\to A$ iff
$e_\Gamma\le[A]_\Gamma$\\
\bottomrule
\end{tabularx}
\end{table}

\subsection{Comparison extensions that preserve nonderivability}

For algebra elements write
\[
 \lambda_c(x)=(c\backslash(xc))\wedge e,
 \qquad \rho_c(x)=((cx)/c)\wedge e.
\]
The two occurrences of a parameter in these individual terms are fixed;
the parameters of the two different terms in PRL remain independent.
Define the comparison formulas
\[
 P(A,B)=((A\vee B)\to A)\wedge t,
 \qquad Q(A,B)=((A\vee B)\to B)\wedge t.
\]

\begin{lemma}[Comparison extension]
\label{lem:comparison-extension}
If $\Gamma\nvdash F$, then for every $A,B$ at least one of
\[
 \Gamma\cup\{P(A,B)\}\nvdash F,
 \qquad
 \Gamma\cup\{Q(A,B)\}\nvdash F
\]
holds.
\end{lemma}
\begin{proof}
Work throughout in the single algebra $L_\Gamma$ and abbreviate its
unit by $e$. PRL and~\eqref{eq:strong-designated} imply
\begin{equation}\label{eq:strong-independent-prl}
 \lambda_c((a\vee b)\backslash a)
 \vee\rho_d((a\vee b)\backslash b)=e.
\end{equation}
Both terms are bounded by $e$, which gives the reverse inequality.
For this semilinearity identity see also~\cite[Section~4, (SL1)]{GLT-Hulls}
and~\cite[Equations~6.8]{BT2003}. Taking $c=d=e$ yields
\begin{equation}\label{eq:strong-comparison-cover}
 (((a\vee b)\backslash a)\wedge e)
 \vee(((a\vee b)\backslash b)\wedge e)=e.
\end{equation}

Two elementary covering facts will control actual derivations with
one additional assumption. First, if $p,z\le e$ and $e\le p\vee z$,
then $p\vee z=e$, so~\eqref{eq:strong-independent-prl} gives
$\lambda_c(p)\vee\rho_d(z)=e$. Taking $d=e$, or interchanging $p,z$
and taking the left parameter to be $e$, gives respectively
\begin{equation}\label{eq:strong-conjugate-cover}
 e\le\lambda_c(p)\vee z,
 \qquad e\le\rho_c(p)\vee z.
\end{equation}
Thus only the active term is conjugated. Second, for $p,q,z,w\le e$,
\begin{equation}\label{eq:strong-product-cover}
 e\le p\vee z,\quad e\le q\vee w
 \quad\Longrightarrow\quad e\le pq\vee(z\vee w).
\end{equation}
Indeed, product preserves finite joins by residuation, and
\[
 e=ee\le(p\vee z)(q\vee w)
   =pq\vee zq\vee pw\vee zw
   \le pq\vee(z\vee w).
\]
The last step uses $zq\le z$, $pw\le w$ and $zw\le z$.
No factors have been permuted.

For a formula $D$, abbreviate $[D]_\Gamma\wedge e$ by $D^-_\Gamma$.
We claim that if
\begin{equation}\label{eq:strong-insert-hypotheses}
 [U]_\Gamma\le e,\qquad z\le e,\qquad
 e\le[U]_\Gamma\vee z,
\end{equation}
then every supplied derivation of $D$ from $\Gamma\cup\{U\}$ yields
\begin{equation}\label{eq:strong-insert-cover}
 e\le D^-_\Gamma\vee z.
\end{equation}
Prove this by induction on that finite derivation. A leaf equal to $U$
uses~\eqref{eq:strong-insert-hypotheses}. An old assumption or an axiom
already has a proof over $\Gamma$, so its class is at least $e$ and
its negative part equals $e$.

At modus ponens, suppose the premises are $D$ and $D\to E$.
Apply~\eqref{eq:strong-product-cover} to their two inductive covers,
using $z$ for both passive terms. Moreover,
\[
 D^-_\Gamma(D\to E)^-_\Gamma
 \le [D]_\Gamma\bigl([D]_\Gamma\backslash[E]_\Gamma\bigr)
 \le[E]_\Gamma.
\]
This product is also at most $e$, because both factors are at most $e$.
Hence it is at most $E^-_\Gamma$, proving the required cover.
The rule $(\mathrm{ADJ}_{\mathrm U})$ uses idempotence of meet with $e$.

For $(\mathrm{PN}_{\to})$, with premise $D$ and parameter $C$, apply the first
inequality in~\eqref{eq:strong-conjugate-cover} to the inductive cover
of $D^-_\Gamma$. Since conjugation is monotone and
$D^-_\Gamma\le[D]_\Gamma$, it follows that
\[
 e\le\lambda_{[C]_\Gamma}(D^-_\Gamma)\vee z
   \le\lambda_{[C]_\Gamma}([D]_\Gamma)\vee z.
\]
The latter active term is precisely the negative part of the rule's
conclusion. For $(\mathrm{PN}_{\rightsquigarrow})$, use the second
inequality and $\rho_{[C]_\Gamma}$ in the same way. This handles all
constructors of the actual derivation and proves the claim.

Suppose now that $[U]_\Gamma,[W]_\Gamma\le e$ and
$e\le[U]_\Gamma\vee[W]_\Gamma$, and that both extensions by $U$ and
by $W$ prove $F$. Applying the claim first with passive term
$[W]_\Gamma$ gives $e\le F^-_\Gamma\vee[W]_\Gamma$. Apply it again
to the proof from $\Gamma\cup\{W\}$, this time with passive term
$F^-_\Gamma$. The result is
\[
 e\le F^-_\Gamma\vee F^-_\Gamma=F^-_\Gamma\le[F]_\Gamma,
\]
so~\eqref{eq:strong-designated} gives $\Gamma\vdash F$.
Taking $U=P(A,B)$ and $W=Q(A,B)$ is permitted by
\eqref{eq:strong-comparison-cover}. Both extensions therefore cannot
prove $F$ when $\Gamma\nvdash F$.
\end{proof}

The proof uses the independent parameters precisely in
\eqref{eq:strong-conjugate-cover}, including the normalization-rule
cases of the induction. It neither assumes that a join equal to $e$
has one summand equal to $e$, nor assumes a chain representation for
the relative algebra.

\subsection{A chain countermodel for an entire theory}

\begin{lemma}[A maximal theory with a chain quotient]
\label{lem:maximal-chain-theory}
If $\Gamma\nvdash F$, there is a theory $\Delta\supseteq\Gamma$ such
that $\Delta\nvdash F$, every proper extension of $\Delta$ proves $F$,
and the order of $L_\Delta$ is total.
\end{lemma}
\begin{proof}
Order the theories containing $\Gamma$ and not proving $F$ by inclusion.
The union of a nonempty inclusion chain is again such a theory.
Otherwise Lemma~\ref{lem:theory-finite-support} extracts a finite list
of assumptions from a proof of $F$ over the union. A finite list in
the union of a nonempty inclusion chain is contained in one member:
successively compare a member containing the next entry with the
member already containing the preceding entries. Rebuilding the proof
over that member contradicts its failure to prove $F$. The empty
chain has the upper bound $\Gamma$. Zorn's lemma therefore gives a
maximal $\Delta$.

For any $A,B$, apply Lemma~\ref{lem:comparison-extension} to $\Delta$.
One of the two comparison extensions still fails to prove $F$; by
maximality the corresponding formula already belongs to $\Delta$.
If it is $P(A,B)$, then
\[
 e_\Delta\le
 \bigl(([A]_\Delta\vee[B]_\Delta)\backslash[A]_\Delta\bigr)
 \wedge e_\Delta.
\]
Residuation implies $[A]_\Delta\vee[B]_\Delta\le[A]_\Delta$, and
therefore $[B]_\Delta\le[A]_\Delta$. If it is $Q(A,B)$, the reverse
comparison follows. Every quotient element has a formula
representative, so the quotient order is total. No additional closure
condition on $\Delta$ is needed: its relative order already uses
derivability over $\Delta$.
\end{proof}

\begin{theorem}[Countable countermodel for arbitrary theories]
\label{thm:theory-chain-countermodel}
In the fixed countable language, if $\Gamma\nvdash F$, then there are
a nontrivial countable bounded unit-cyclic residuated chain $C$ and a
single valuation $v$ such that
\[
 e_C\le v(A)\quad(A\in\Gamma),\qquad v(F)<e_C.
\]
When the language contains $f$, the same statement includes an
interpretation $f_C$ of that independent constant.
\end{theorem}
\begin{proof}
Take $\Delta$ from Lemma~\ref{lem:maximal-chain-theory} and put
$C=L_\Delta$. Lemma~\ref{lem:relative-lindenbaum} supplies its full
bounded residuated structure and unit cyclicity; maximality supplies
the total order. Its carrier is countable, since it is a quotient of
the countable set of finite formulas. By
\eqref{eq:strong-designated}, $\Delta\nvdash F$ means
$e_C\nleq[F]_\Delta$, hence $[F]_\Delta<e_C$. In particular, $C$ is
nontrivial.

Use the canonical valuation $v(p)=[p]_\Delta$, and set
$f_C=[f]_\Delta$ if needed. Formula induction gives
$v(A)=[A]_\Delta$ for every $A$. Every premise in $\Gamma$ belongs
to $\Delta$ and thus has a hypothesis proof there. Equation
\eqref{eq:strong-designated} makes all these premises designated
under this one valuation, whereas $F$ is strictly below the unit.
\end{proof}

Together with the soundness of the Hilbert rules, this establishes
chain completeness for arbitrary theories of the explicitly specified
calculus. Its construction uses finite support of proof trees, not
compactness of chain or standard semantics.

\subsection{Full standard embeddings and strong completeness}

The remaining algebraic input must apply to chains that need not be
dense. The dense-chain completion used for theoremhood earlier in the
paper does not by itself provide such an input.

\begin{theorem}[Full bounded standard embedding {\cite[Corollary~10.5]{Wang2026}}]
\label{thm:countable-standard-embedding}
Every nontrivial countable bounded unit-cyclic residuated chain $C$
has an embedding $j:C\hookrightarrow S$ into a standard algebra,
preserving
\[
 \wedge,\ \vee,\ \cdot,\ \backslash,\ /,\ e,\ \bot,\ \top.
\]
No density assumption on $C$ is required.
\end{theorem}
\begin{proof}
This is the bounded representation theorem
of~\cite[Corollary~10.5]{Wang2026}. Its construction first uses a
cyclic involutive extension preserving the original residuated
signature, applies the full standard representation, and then
restores the original bounds by an interval restriction and rescaling.
The resulting map retains both actual residuals. We use the statement
of the specified version arXiv:2611.19940v1: its hypotheses are precisely
nontriviality, countability, boundedness and unit cyclicity. In particular,
neither density nor involutivity is required of the original $C$.
The representation construction itself is an external input to this
proof; the companion implementation does not replace that citation.
\end{proof}

The independently interpreted constant $f$ imposes no further
requirement on this theorem: after choosing $j$, define its value in
the target to be $j(f_C)$. This does not identify $f_C$ with the
dualizing constant of an auxiliary involutive extension.

\begin{theorem}[Strong standard completeness]
\label{thm:strong-standard}
For every theory $\Gamma$ and formula $F$ in the fixed countable
language of $\mathrm H_{\mathrm{sl}}^{*}$,
\[
 \Gamma\vdash_{\mathrm H_{\mathrm{sl}}^{*}}F
 \quad\Longleftrightarrow\quad
 \Gamma\models_{\mathsf{Std}}F.
\]
The theory may be infinite. The result holds in the full language
with the independent constant $f$, with all its interpretations
included in the standard semantics.
\end{theorem}
\begin{proof}
For soundness, fix a standard algebra and a valuation satisfying
$\Gamma$. A1--A18 are sound in bounded residuated lattices and WCM
is sound by weak exchange. For PRL, totality gives $a\vee b=a$ or
$a\vee b=b$. In the first case $e\le(a\vee b)\backslash a$, so its
left conjugate, truncated at $e$, equals $e$ for every parameter; in
the second case the right conjugate equals $e$ for every parameter.
Thus independent PRL is sound.

The assumption leaves are designated by the choice of valuation.
Modus ponens preserves designation because $e\le a\backslash b$
implies $a\le b$. Conjoining a designated value with $e$ preserves
designation. Finally, $e\le a$ gives $c=ec\le ac$ and
$c=ce\le ca$, so residuation gives
$e\le c\backslash(ac)$ and $e\le(ca)/c$. These are the two
normalization rules. Induction on a finite derivation proves the
forward implication.

For the converse, suppose $\Gamma\nvdash F$. Obtain $C,v$ from
Theorem~\ref{thm:theory-chain-countermodel} and a full bounded
embedding $j:C\hookrightarrow S$ from
Theorem~\ref{thm:countable-standard-embedding}. Set
$v_S(p)=j(v(p))$ for every ordinary variable, and $f_S=j(f_C)$ if
$f$ is present. Preservation of the entire signature gives, by
induction on formulas,
\begin{equation}\label{eq:strong-evaluation-transfer}
 v_S(A)=j(v(A))\qquad\text{for every formula }A.
\end{equation}
Since $j(e_C)=e_S$ and $j$ preserves order, every member of $\Gamma$
remains designated. An order embedding also preserves strict order,
so $v_S(F)<e_S$. This one standard algebra and one valuation refute
$\Gamma\models_{\mathsf{Std}}F$.

The same embedding is used for all premises, including an infinite
set of them. Formula induction requires preservation of the finitary
operations and constants, not preservation of infinite joins or
meets. This establishes the asserted arbitrary-theory result.
\end{proof}

\begin{corollary}[Compactness of standard consequence]
\label{cor:standard-compactness}
For every $\Gamma,F$ in this countable language,
\[
 \Gamma\models_{\mathsf{Std}}F
 \quad\Longleftrightarrow\quad
 \text{there exist $A_1,\ldots,A_n\in\Gamma$ such that }
 \{A_1,\ldots,A_n\}\models_{\mathsf{Std}}F,
\]
where the finite list may be empty.
\end{corollary}
\begin{proof}
Apply Theorem~\ref{thm:strong-standard}, then
Lemma~\ref{lem:theory-finite-support}, and then soundness to the
finite supported derivation. Conversely, satisfying all members of
$\Gamma$ satisfies the displayed finite set, so the implication
follows directly from the definition of semantic consequence.
\end{proof}

Thus semantic compactness is a consequence of strong completeness,
and was not an input to the maximal-theory construction. The result
of this section concerns Hilbert consequence from specified premises.
Theorem~\ref{thm:hilbert-explicit} remains a correspondence for
formula theorems and does not, without a further relative translation,
identify this consequence relation with an arbitrary-assumption
hypersequent version of $\G$. No full cut-elimination conclusion is
used or asserted.

\section{Scope of the correction}\label{sec:scope}

Table~\ref{tab:corrections} records the relation between the published
argument and its replacement. Original item numbers refer throughout to
the published version of~\cite{Wang2019}, rather than to a preprint or
to the separate paper on involutive uninorm logic.

\begin{table}[htbp]
\centering\small
\caption{Changes to the proof and the conclusions retained.}
\label{tab:corrections}
\begin{tabularx}{\linewidth}{@{}>{\raggedright\arraybackslash}p{.22\linewidth}>{\raggedright\arraybackslash}X>{\raggedright\arraybackslash}X@{}}
\toprule
Original item & Required change & Result established here\\
\midrule
Definition 5(iv) & One occurrence matching and one renaming on both sides & Producer--consumer incidence preserved\\[3pt]
Theorem 2 & Supply the omitted normalization proof & Theorem~\ref{uc:normalization}\\[3pt]
Construction 3 and Lemma 5 & Replace the density image and auxiliary contraction & Theorem~\ref{thm:cartesian-local}\\[3pt]
Lemma 7 & Replace the global separation construction & Theorem~\ref{olg:source-free}\\[3pt]
Theorem 4 & Use the new normalization and separation arguments & Theorem~\ref{thm:strong-density}\\[3pt]
Theorem 5 & Supply the semantic argument and retain theorem-level scope & Corollary~\ref{cor:fss-weak-standard}\\
\bottomrule
\end{tabularx}
\end{table}

The finite countermodel concerns the compatibility of a labelled
contraction with an ordered density transformation. It does not refute
the original failure of full cut elimination. The normalization proved
here retains cuts on $t$, exactly where the producer-side interaction
with WCM prevents the unrestricted cut-free argument.

The Cartesian construction also changes the assertion required of the
global algorithm. It is enough to produce a nonempty closed theorem
containing only terminal components, using contractions already justified
by the local invariant. Uniqueness of intermediate fragments and the
original particular schedule of constrained contractions are not needed.
The preparation and separation proofs in the appendices establish that
weaker, sufficient interface with explicit occurrence histories.

Combining density admissibility with the semantic construction yields
weak standard completeness of $\G$. Appendix~\ref{sec:hilbert-explicit}
gives the Hilbert axioms for $\mathrm H_{\mathrm{sl}}^{*}$ and proves
that it has the same formula theorems as $\G$. This supplies the proof
of the correspondence indicated by circled link~1 in Figure~4
of~\cite[p.~12, proof of Theorem~5]{Wang2019}.

The semantic construction in Section~\ref{sec:fss-semantics}
starts from a countable dense countermodel for the density-augmented calculus.
Theorem~\ref{thm:strong-standard} concerns arbitrary theories of the explicit
$\mathrm H_{\mathrm{sl}}^{*}$ in a countable language, with all premises
evaluated at the same valuation. Its direct relative chain-countermodel
argument is combined with the external bounded representation in
Theorem~\ref{thm:countable-standard-embedding}. A corresponding strong
standard completeness result was already stated
in~\cite[Corollary~10.7]{Wang2026}.
It does not assert a corresponding arbitrary-theory translation for $\G$,
completeness with respect to one fixed standard algebra, or full cut
elimination.

\paragraph{Proof-development and verification methods.}
OpenAI ChatGPT was used interactively to propose candidate counterexamples,
develop and critically examine proof transformations, compare proof
dependencies, and develop and revise the accompanying Lean~4 code.
Separate model-assisted checks examined local rule cases, side conditions,
and induction interfaces. These exchanges supplied candidates for
mathematical assessment; agreement between model outputs was not treated
as a proof criterion. The mathematical justifications are the explicit
arguments presented in this manuscript. For the formalized results,
Lean's kernel checks the elaborated proof terms. The archived release
\emph{GpsUL-1.0.0}~\cite{Wang2026GpsUL} records the toolchain and
dependencies and provides build instructions, axiom-audit scripts, and
the correspondence between manuscript statements and Lean declarations.
Kernel checking establishes the formal statements; correspondence with
the mathematical prose remains a separate review obligation.

\appendix
\renewcommand{\thesection}{\Alph{section}}
\setcounter{table}{0}
\makeatletter
\def\appendixprefix@section{Appendix~}
\def\@seccntformat#1{%
  \ifcsname appendixprefix@#1\endcsname
    \csname appendixprefix@#1\endcsname\fi
  \csname the#1\endcsname.\quad}
\makeatother

\section{Normalization to unit cuts}
\label{uc:section}

The density argument requires a proof in which the only cut formula is
the monoidal unit. We first establish this normalization independently
of density elimination. Antecedents remain ordered words throughout.
The producer and receiver of a cut are, respectively, its premises
\[
 G_0\mid\Pi\seq A,
 \qquad H_0\mid\Gamma,A,\Delta\seq B,
\]
and its conclusion is
$G_0\mid H_0\mid\Gamma,\Pi,\Delta\seq B$.
The displayed antecedent occurrence of $A$ is distinguished from every
other occurrence of the same formula.

\begin{theorem}[Unit-cut normalization]\label{uc:normalization}
For each finite $\G$ proof of a hypersequent $H$ there is a finite
$\Gtwo$ proof of $H$. Consequently these calculi have the same theorems.
\end{theorem}

Write $|A|$ for the number of binary connectives in $A$ and
$\operatorname{Adm}(A)$ for admissibility in $\Gtwo$ of the displayed
cut on $A$, for arbitrary contexts and premise proofs. Put
\[
 \operatorname{LowerCuts}(A)\quad\Longleftrightarrow\quad
 \forall C\,\bigl(|C|<|A|\ \Longrightarrow\ \operatorname{Adm}(C)\bigr).
\]
This hypothesis ranges over all formulas of smaller complexity, including
those that do not occur in $A$. Fix $A\ne t$ and assume
$\operatorname{LowerCuts}(A)$. Prove the left substitution lemma by
induction on the receiver proof and then the right substitution lemma
by induction on the producer proof. Their principal reductions below
invoke $\operatorname{LowerCuts}(A)$ only on immediate subformulas of $A$.

A mark selects a whole formula occurrence, not every occurrence with
the same written form. Backward tracing duplicates its pattern in
shared lattice contexts and in the two premise copies of external
contraction. Newly inserted formulas are not marked. We also use the
following elementary bookkeeping: a prescribed external background
can be added by EW, and several complete copies of that background
can be reduced to one by componentwise EC. This preserves any
multiplicity prescribed in the background itself; it never projects
an arbitrary component from a hypersequent.

\begin{lemma}[Left substitution]\label{uc:left}
Fix $A\ne t$ and assume $\operatorname{LowerCuts}(A)$.
Let $P$ prove $G_0\mid\Pi\seq A$ in $\Gtwo$, ending in a principal
right introduction of $A$, or in the initial sequent $\Pi\seq\top$
when $A=\top$. If $Q$ proves $K$ in $\Gtwo$, mark any set of antecedent
occurrences of $A$ in $K$. Replacing every marked occurrence by the
word $\Pi$ yields a hypersequent $K[\Pi]$ such that
\[
                 \Gtwo\vdash G_0\mid K[\Pi].
\]
The number of marks may be zero.
\end{lemma}

\begin{proof}
Induct on the height of $Q$, for the fixed principal producer $P$ of
arbitrary height. With no marks, add $G_0$ by EW. For an initial
identity, a mark is possible only in $A\seq A$, whose transform is
provided by $P$. A transformed $\top$-right initial sequent is again
initial. In a $\bot$-left initial sequent the distinguished $\bot$
cannot be marked, because there is no principal right producer of the
kind stipulated for $A=\bot$; hence $\bot$ remains and the transform
is initial. The unit initial sequent has no antecedent mark. Add
$G_0$ by EW where necessary.

For EW, transform its premise and add the transformed new component.
For EC, mark both matching premise copies and apply the induction
hypothesis once; their transformed copies are identical and contract.
For COM, trace the marks to their positions in its two premises and
replay communication after the induction hypotheses: replacement
commutes with concatenation, so each selected middle interval remains
a consecutive word. For WCM, transform each of its two word blocks
and perform the same rotation; its succedent stays $t$.

For WCT, neither displayed cut occurrence of $t$ is marked. The
induction hypotheses give transformed active premises
\[
 \Theta',t,\Lambda'\seq D,
 \qquad \Sigma'\seq t.
\]
Their original unit cut yields
$\Theta',\Sigma',\Lambda'\seq D$.
This includes marks in the receiver contexts, in the producer word
$\Sigma$, or in both. Full external backgrounds accompany these
premises, and duplicated copies of $G_0$ are contracted.

Every logical rule whose principal antecedent formula is unmarked
is replayed after tracing only the marked context occurrences.
Its principal formulas and immediate subformulas remain unmarked.
For residual right introduction, in particular, the discharged
boundary formula is unmarked even if it is syntactically $A$, so it
stays at the required end of the antecedent. Residual left introduction
inserts its transformed auxiliary word at its original position;
product right introduction concatenates its transformed words.
The common contexts of $\wedge R$ and $\vee L$ inherit identical
patterns in both premises and therefore still match. This also covers
product left, meet left, join right, and unit left. Contract the
external copies of $G_0$ supplied by different premises.

It remains that the last rule of $Q$ introduces a marked antecedent
$A$. First use the height induction in its premises to replace all
\emph{other} marks. The newly exposed proper subformulas of $A$ are
unmarked. Primes below denote transformed contexts; external
backgrounds are suppressed only in the displays and are carried by
every cut.

If $A=U\cdot V$, the producer premises are
$G_1\mid\Pi_1\seq U$ and $G_2\mid\Pi_2\seq V$, with
$\Pi=\Pi_1,\Pi_2$ and $G_0=G_1\mid G_2$.
Apply the admissible cuts on $U$ and $V$ successively to the receiver
premise $\Gamma',U,V,\Delta'\seq B$, obtaining
\[
                \Gamma',\Pi_1,\Pi_2,\Delta'\seq B.
\]
For $A=U\backslash V$, the producer premise is
$G_0\mid U,\Pi\seq V$. The transformed receiver premises have
active sequents $\Sigma'\seq U$ and
$\Gamma',V,\Delta'\seq B$. Cutting first on $U$ and then on $V$
gives $\Gamma',\Sigma',\Pi,\Delta'\seq B$, exactly the replacement
in $\Gamma',\Sigma',U\backslash V,\Delta'\seq B$.
For $A=V/U$, start instead with $G_0\mid\Pi,U\seq V$; the same two
lower cuts give $\Gamma',\Pi,\Sigma',\Delta'\seq B$, the replacement
in $\Gamma',V/U,\Sigma',\Delta'\seq B$. No exchange is used.

For $A=U\wedge V$, use the producer premise for the conjunct selected
by the receiver's meet-left rule, followed by its lower cut. The
induction hypothesis on the receiver premise already included the
whole producer background $G_0$, including that of the unused producer
premise. Contract any duplicate background supplied by the used premise.

For $A=U\vee V$, use the receiver premise corresponding to the
disjunct selected by the producer's join-right rule. After transforming
its other marks, cut on that proper subformula. Add by EW the
syntactically transformed external background of the unused receiver
premise. This recovers precisely the background of the original
join-left conclusion; its unused active component is not required.
In particular, no projection from that unused premise is made.

There is no principal negative case for $\top$; the cases of variables
and $\bot$ have no stipulated principal producer, and $t$ is excluded.
Thus every rule is covered. Same-formula recursive calls used strict
subproofs of $Q$, and every displayed cut used a proper subformula.
\end{proof}

\begin{lemma}[Right substitution]\label{uc:right}
Fix $A\ne t$ and assume $\operatorname{LowerCuts}(A)$.
Fix a $\Gtwo$ proof $Q$ of
$H_0\mid\Lambda,A,\Theta\seq B$ with the displayed antecedent
occurrence designated. If a $\Gtwo$ proof $P$ ends in
\[
 G_0\mid\Pi_1\seq A\mid\cdots\mid\Pi_m\seq A,
\]
where precisely the displayed succedent occurrences are marked, then
\begin{equation}\label{uc:right-result}
 \Gtwo\vdash G_0\mid H_0\mid
 \Lambda,\Pi_1,\Theta\seq B\mid\cdots\mid
 \Lambda,\Pi_m,\Theta\seq B.
\end{equation}
Here $m$ may be zero.
\end{lemma}

\begin{proof}
Induct on the height of $P$, with $Q$ fixed. With no marks, add $H_0$
by EW. Call the replacement
$\Sigma\seq A\mapsto\Lambda,\Sigma,\Theta\seq B$ its receiver
transform. A marked initial identity is handled by $Q$ itself.
A marked $\bot$-left initial component still contains $\bot$ after
its receiver transform and is therefore initial; add $H_0$ by EW.
A marked $\top$-right initial component is handled by
Lemma~\ref{uc:left}, already established at this same formula.
The initial $\seq t$ cannot be marked.

At EW, apply the induction hypothesis and add the new component,
transformed if marked. At EC, a marked conclusion marks both contracted
premise copies. A single induction-hypothesis call transforms both
copies identically, and EC applies. An unmarked contraction is replayed.
This simultaneous treatment avoids sequential cuts whose intermediate
proof heights would not be bounded by the original height.

At COM, trace each marked output succedent to its premise and apply
the induction hypotheses. Wrap the antecedent of each marked active
component with $\Lambda,\Theta$, and keep an unmarked active component
unchanged. Exchange exactly the original middle intervals. For example,
when both active succedents are marked, communication transforms
\[
 \Lambda,\Gamma_1,\Pi_1,\Delta_1,\Theta\seq B,
 \qquad
 \Lambda,\Gamma_2,\Pi_2,\Delta_2,\Theta\seq B
\]
into their counterparts with $\Pi_1,\Pi_2$ exchanged, exactly as
required. With one or no mark, omit the wrappers for the unmarked
components and retain their succedents. Contract duplicate $H_0$.

The active succedent of WCM is $t$ and is unmarked, since $A\ne t$.
Only its background is transformed, and WCM is replayed. At WCT, the
minor succedent $t$ is likewise unmarked. If the major succedent is
marked, its transformed active sequent is
\[
                 \Lambda,\Gamma,t,\Delta,\Theta\seq B.
\]
Cutting its displayed $t$ with the unchanged minor active sequent
$\Sigma\seq t$ yields
$\Lambda,\Gamma,\Sigma,\Delta,\Theta\seq B$.
All marked backgrounds in both premises are handled by their induction
hypotheses. If the major succedent is unmarked, replay the original
unit cut on unchanged active components. Remove duplicate $H_0$.

A logical rule with an unmarked active conclusion transforms only
external marks and is replayed unchanged. If its active succedent is
marked but is not logically introduced, the rule is a left rule or
unit left, and is replayed within the added outer words
$\Lambda,\Theta$. Product left, meet left, and unit left retain their
displayed positions. In residual left introduction only the major
succedent inherits the active mark; the auxiliary active succedent
remains unmarked, even when its formula happens to equal $A$.
The auxiliary antecedent is inserted at its original position in
the enlarged major context. At join left both active premise
succedents inherit the mark and their enlarged common contexts match.
These cases cover all left rules; remove any duplicate $H_0$ backgrounds.

Finally, suppose the last rule introduces a marked succedent $A$ on
the right. Apply the height induction to its premises, marking only
the ancestors of the \emph{other} marked succedents in the external
background. The active premise succedents are unmarked. Replay the
same right introduction; its antecedents, including the common word
for meet right, are unchanged. This yields a proof
\[
                         P': F\mid\Pi\seq A,
\]
where $F$ contains the required unchanged backgrounds, the receiver
transforms of the other marked components, and copies of $H_0$.
Leave these copies uncontracted for now, so that $P'$ ends
\emph{literally} in the principal right introduction.
Apply the already proved Lemma~\ref{uc:left} to $P'$ and the one
designated antecedent occurrence in the original $Q$. Its conclusion is
\[
                     F\mid H_0\mid\Lambda,\Pi,\Theta\seq B.
\]
Contract the duplicate backgrounds to obtain
\eqref{uc:right-result}. Crucially, no application of the present
lemma is made to $P'$: its only recursive calls were to strict
premises of the original $P$.
\end{proof}

\begin{proof}[Proof of Theorem~\ref{uc:normalization}]
Prove $\operatorname{Adm}(A)$ for every formula $A$ by strong induction
on $|A|$. At the current complexity, the induction hypothesis supplies
$\operatorname{LowerCuts}(A)$. For $A=t$, cut is WCT. For $A\ne t$,
instantiate Lemma~\ref{uc:right} with this hypothesis and $m=1$.
Its conclusion is the required cut on $A$ with its two given premise
proofs. The lemma uses Lemma~\ref{uc:left} at $A$, which was proved
first without invoking $\operatorname{Adm}(A)$. For variables,
$\bot$ and $\top$, the complexity is zero and
$\operatorname{LowerCuts}(A)$ is vacuous.

Now normalize an arbitrary finite $\G$ proof by induction on its proof
tree. Normalize the premises of its last inference. Replay any noncut
inference, retain a cut on $t$ as WCT, and replace a nonunit cut by the
admissible construction just proved. The end-hypersequent is unchanged.
The reverse inclusion of theorems holds because WCT is a special cut.

For the explicit proof transformation, carry a cut-replacement procedure
for each smaller formula complexity through the outer induction. The
two substitution procedures use structural recursion on the original
receiver and producer proofs, respectively. The right procedure may
invoke the already defined left procedure at the same formula.
The final recursion follows the original finite proof tree.
\end{proof}

The archived declaration \texttt{unit\_cut\_normalization} returns
\texttt{KProvable}~$H$, defined as \texttt{Nonempty} of the type of
unit-cut proof trees of $H$. Its checked conclusion is the finite-proof
existence assertion above. The effective transformation is specified
by the preceding recursions on annotated finite trees.

Theorem~1 of~\cite{Wang2019} exhibits a theorem of $\G$ with no cut-free
proof. By Theorem~\ref{uc:normalization} it has a proof using only unit
cuts, so every such proof of that theorem contains at least one unit cut.
The preprocessing in Appendix~\ref{sec:ordered-preparation} starts from
the unit-cut proof supplied by Theorem~\ref{uc:normalization}.

\section{Ordered preprocessing and occurrence histories}\label{sec:ordered-preparation}
\subsection{From unit cuts to a linear labelled proof}
\label{op:section}

Let $\pi$ be a finite proof, containing only cuts on the unit $t$, of
\[
 G_0=G'\mid\{\Gamma_i,p,\Delta_i\Rightarrow A_i\}_{i=1}^{n}
          \mid\{\Pi_j\Rightarrow p\}_{j=1}^{m},
\]
where the displayed occurrences are the only occurrences of $p$ in $G_0$.
In particular, every component of $G_0$ has at most one such occurrence.
For an arbitrary $\G$ proof of this premise, Theorem~\ref{uc:normalization}
supplies such a $\pi$. Antecedents are ordered words, and equal components
at different positions remain distinct occurrences. Use the generalized lattice rules
\eqref{eq:cart-generalized-meet} and~\eqref{eq:cart-generalized-join} of
Definition~\ref{def:cart-auxiliary}. Their equal-context instances followed
by ordinary external contraction give the corresponding ordinary lattice
rules. These generalized rules remain primitive inside the labelled
constructions; their unlabelled derivations are not inlined there.

\begin{proposition}[Preparation from an actual unit-cut proof]
\label{op:preparation}
From $\pi$ one can construct a finite proof $\tau$ using the generalized
lattice rules, WCM and WCT, but no external weakening or external
contraction.  Its nonempty closed conclusion has the form $R=G\mid G^*$.
It carries the following occurrence histories, which will control all later
copies and extractions.
\begin{enumerate}
\item Every source is the position of a designated, actual premise-focus
occurrence of a noncontraction, noninert rule.  It is the terminal focus of
a component thread containing at least one contraction in the precursor
proof.  An identity leaf cannot be declared a source merely by attaching a
tag to it.
\item A residual of source $s$ is an occurrence frozen at an omitted
contraction of that particular thread.  It is distinct from the live focus
at $s$.  It remains background in $\tau$ from its freeze position onwards.
All residuals of source $s$ are already frozen by the time the live focus
at $s$ is used.
\item The thread ends at that live focus. Each principal component
produced by the rule starts a new thread. Each background component
transmitted unchanged continues its own thread. Consequently no continuation
below that focus can newly freeze a residual with the same source $s$.
\item If $b(R_0)$ is the freeze position of a residual $R_0$ and $s(R_0)$
is its source, then $s(R_0)\le b(R_0)$ in the root-minimal tree order.
If a complete auxiliary subtree contains $R_0$ at its root but does not
contain $s(R_0)$, then $s(R_0)$ lies on the continuation from that subtree's
root to the global root.  Copies retain the map to this original history.
\end{enumerate}
Every component of $G$ has a terminal reference to a component occurrence
of $G_0$, and $G$ is nonempty.  These terminal components have at most one marked
slot.  A component of $G^*$ has the ordered slot skeleton of its source
focus.  The source set is finite.
\end{proposition}

\begin{proof}
We give the construction and all history assignments.

\emph{First, prune external weakening.}
Induct from the leaves of $\pi$ to its root.  At each node $K$ retain a
proof of a nonempty component submultiset $K^-\subseteq K$, with its
occurrence inclusion.  Keep initial sequents.  At $B/(B\mid S)$, keep
only the already constructed proof of $B^-$.  At an external contraction,
contract precisely the retained designated copies if at least two remain;
if one remains, map it to the designated survivor without inference; if
none remains, retain the nonempty surviving background.  This uses the
ordered equality of the contracted sequents, not a permutation of their
antecedents.

At a unary noncontraction rule, replay it if its focus survives; otherwise
return its retained background.  At a binary rule, replay it if both
focuses survive.  If a focus in one premise does not survive, return that
premise's retained nonempty background and discard the other proof.
This background is present unchanged in the original conclusion.
The same argument applies to WCT: when both focuses survive it is the
same unit-cut instance, and otherwise the focus-free premise supplies
unchanged background.  It also applies to ordered COM and the ordinary
lattice rules, whose complete focus sequents retain the required common
contexts whenever both focuses survive.  The resulting proof $\pi^-$
has no external weakening, still has only unit cuts, and ends in a
nonempty $G_0^-\subseteq G_0$.

\emph{Second, linearize the lattice contexts.}
Replace each ordinary $\wedge_r$ or $\vee_l$ by its generalized version
followed by contraction of the two identical output components.  This
leaves the original premises and conclusion unchanged and introduces no
cut.  Call the result $\widehat\pi$.  The generalized rules now transmit
the context occurrences of their two premises separately.

Every formula occurrence in $\widehat\pi$ containing $p$ is a bare $p$.
To see this first in $\pi^-$, follow a formula occurrence towards the
root: an ordinary logical rule either retains it or makes it a subformula
of its successor; contraction identifies it with its matching survivor;
COM and WCM only redistribute occurrences.  A unit cut can delete an
occurrence of $t$, but cannot delete a formula containing $p$.  A proper
occurrence of $p$ would therefore have a proper occurrence at the root,
contrary to the form of $G_0^-$.  The lattice replacement adds only copies
of formula shapes already present.  Thus no occurrence of $p$ can be a
logical principal formula in a retained inference.

\emph{Third, define the precursor threads before omitting contractions.}
Assign distinct identities to all component occurrences of
$\widehat\pi$.  Join each background occurrence to its unchanged output
occurrence.  At a contraction, join every designated contracted occurrence
to the designated survivor. For any other noninert rule, each premise
focus is the terminal occurrence of its thread. Each newly produced principal
component starts a new thread, and the background edges continue the
threads of the unchanged components. Root components are also terminals.
An inert position marker preserves its thread.

Each component occurrence has a unique terminal under these edges, since
it has at most one rootward successor and the proof is a finite tree.
A thread class consists of all occurrences having the same terminal.
All components in a thread class have the same ordered sequent shape:
its edges are background transmission or contraction of identical ordered
sequents.  At contraction match equal formula occurrences by their actual
positions in the antecedent word and the unique succedent.  Repeated bare
$p$ occurrences are thereby matched by positions, not arbitrarily
permuted.

For each thread class containing a contraction and ending at an actual
noncontraction premise focus, declare that focus to be one source $s$.
Its common sequent shape, with bare $p$ positions marked, is its skeleton
$\sigma_s$.  Different focus occurrences define different sources even if
their sequent shapes agree.  A thread ending at the root has a terminal
reference instead.  Insert an inert marker on the appropriate premise
edge to distinguish source positions in different premises of one rule.
This refines the finite proof tree and gives the root-minimal order on
sources. An identity that immediately feeds a noncontraction rule has
no contraction in that thread and therefore does not define a source.
More generally, no thread without a contraction defines a source.

\emph{Fourth, omit contractions by freezing their discarded occurrences.}
Process $\widehat\pi$ from the leaves towards its root, always retaining the
occurrences designated by the original proof.  At each contraction,
omit the inference, use its designated survivor as the representative of
the original output, and mark every other designated contracted
occurrence as frozen.  Already frozen occurrences remain distinct
background components.  At each other inference, replay the original
rule on its original focus occurrences and carry all frozen components
as extra external background.  A binary rule takes the multiset union of
the extra backgrounds.  Every listed rule, including WCM and WCT, allows
this extension of its external contexts.  Hence no new weakening
inference or leaf is required.

Frozen occurrences are never selected as original focuses.  A frozen
occurrence retains the thread class of the precursor occurrence it
represents.  If the class ends at a source $s$, it becomes a residual of
source $s$.  If the class ends at a root component, it is given that
terminal reference.  Original surviving root components also keep their
terminal references.  These are the two exhaustive classes of root
occurrences: a class containing a frozen occurrence contains a
contraction by construction.

The source focus remains an actual designated live occurrence in this
new proof.  It is distinct from all frozen copies bearing its source tag.
The precursor thread stops at that focus, so every omitted contraction
of the class occurs before that stop. The principal components produced
at the focus belong to new threads. A residual of the old source can only
continue as a previously frozen background occurrence.  The construction thereby supplies the
first three history conditions, rather than imposing them on arbitrary
tags after the proof has been built.

For the fourth condition, the path from a freeze position $b(R_0)$ to its
thread terminal contains the source position $s(R_0)$, so
$s(R_0)\le b(R_0)$.  Suppose a complete subtree rooted at $v$ contains
$R_0$ at its root.  Its freeze position is in that subtree.  The unique
path from that position to the global root either meets its source
before leaving the subtree, or leaves through $v$ and meets the source
on the rootward continuation from $v$.  These alternatives prove the
claimed import dichotomy.  In particular, an imported residual whose
source is inside the imported subtree can already have passed its source;
it must not be confused with a residual whose designated source is still
in the future.  Only the latter has its source on the receiving
continuation.

A continuation retained from an actual selected focus follows the original
rootward path.  If a thread created on that continuation subsequently
produces a residual, its terminal source is on that path.  Threads with
no such residual need not define sources.  This is the precise source-routing property used
by extraction; it does not assert that an auxiliary residual with a past
source first formed on the receiving path.  Whole-proof copies retain the
map to this precursor tree, its thread terminals and its positions.

\emph{Fifth, decorate slots and label only paired identity origins.}
Mark every original bare $p$ position, including those in frozen
occurrences, as a slot.  After contractions are omitted, each existing
slot has exactly one successor at every retained inference, on the same
side.  No logical connective contains a slot.  WCM moves its specified
two ordered blocks without changing their internal order; COM exchanges
its specified blocks; WCT deletes only its ordinary $t$ positions.
The generalized lattice rules transmit their two context families
separately.  This proves the required same-side linearity directly.

A bare slot originates either in a paired identity $p\Rightarrow p$ or
in a nonidentity initial context.  The latter can occur in $\bot_l$ or
$\top_r$; there is no external weakening left.  Replace such a
nonidentity-born left slot and all its successors by $\top$, and such a
right slot and all its successors by $\bot$.  Its initial sequent remains
an instance of the same initial rule.  All later rules remain legitimate
because this slot is only a context occurrence.  Preserve the slot mark:
an extreme constant filling a slot is distinguished from an original,
unmarked extreme constant.

Assign a globally fresh label $a$ to each remaining identity leaf
$p\Rightarrow p$, obtaining $p_a\Rightarrow p_a$, and propagate it along
the two unique successor paths.  Each live label has exactly one left
and one right occurrence at every later node on the path from its leaf.
At a binary rule its two premise subtrees have disjoint identity leaves,
and hence disjoint label sets.  It follows by induction that every node
is closed and every binary rule has label-disjoint closed premises.
All nonidentity initial sequents are free of live labels.

Let $G$ be all root occurrences with terminal references, including the
surviving $G_0^-$, and let $G^*$ be the source residuals.  The first is
nonempty.  Its skeletons are those of components of $G_0^-$ and hence
have at most one slot.  Every residual has the ordered skeleton of its
source focus.  Decorations of different occurrences of the same
skeleton need not coincide: only their marked positions and source
history agree.

Finally, a full auxiliary-proof copy is made by a single fresh bijection
on both occurrences of each live label, retaining its constant fillings,
slot marks and precursor-history map.  It is again a proof with closed
nodes and disjoint binary premises.  This realizes the complete-import
convention of the prepared interface.  The construction supplies actual
history data, not an arbitrary declaration that a component is residual.
\end{proof}

\begin{remark}[The history condition]
\label{op:notarbitrary}
The input to ordered separation includes the original unit-cut proof of
$G_0$ and the preparation obtained from it by Proposition~\ref{op:preparation}.
The bound of one bare $p$ per component is imposed on every component
of $G_0$, including components subsequently removed in pruning.
The source tags are determined by the precursor contraction threads.
The actual noncontraction source focus, its distinct
frozen copies, the termination of their precursor thread, and the
nonempty terminal range are part of the input certificate.  The data are obtained entirely from the precursor proof, before source
separation is performed.
\end{remark}
\providecommand{\Om}{\mathsf A_{\mathrm{ord}}}
\providecommand{\ECs}{\mathrm{EC}_{\mathrm{coh}}}
\providecommand{\Par}{\operatorname{Par}}
\providecommand{\supp}{\operatorname{supp}}
\providecommand{\lca}{\operatorname{lca}}

\subsection{The prepared ordered proof interface}

Antecedents are finite ordered words and succedents are single formulas.
Only the order of hypersequent components is immaterial. A slot always means
one specified occurrence in such a word or succedent, and is never an
unlocated occurrence in a multiset. A labelled hypersequent is closed when every live label occurs exactly
once in an antecedent and exactly once in a succedent. The auxiliary
calculus $\Om$ is obtained from the ordered calculus $\G$ as follows:
external weakening is omitted; the lattice rules are replaced by
$(\wedge_{rw})$ and $(\vee_{lw})$; only unit cuts are retained; and every
binary inference has closed premises with disjoint label sets. The
characteristic variable occurs only as a bare sequent formula, paired
identities are $p_a\Rightarrow p_a$, and other initial sequents contain no
live label. Ordinary external contraction is replaced by
\[
 \frac{H\mid C\mid C^{\sigma}}{H\mid C}\ (\ECs).
\]
Here $C$ and $C^{\sigma}$ are disjoint closed blocks, and one bijection
$\sigma$ renames each label simultaneously at its antecedent and succedent
occurrence. The correspondence must preserve the ordered antecedents, ordinary
formulas, component occurrences, slot positions, and all source and protection
data introduced below. The bijection acts on an entire block; it is not chosen
independently on the two sides. This is the coherent contraction rule. Its density compatibility follows
from Theorem~\ref{thm:cartesian-local}: make identical cycle choices in the
two blocks and contract the resulting duplicate sequents in $\G$.

\begin{definition}[Prepared ordered proof]\label{ols:prepared}
Start with a finite unit-cut proof of an end hypersequent $G_0$ in which
$p$ occurs only as a bare formula and each component contains at most one
occurrence of $p$, counting both sides together. Apply the preprocessing
of Proposition~\ref{op:preparation}. Retain its contraction-free labelled
proof tree $\tau$, its nonempty closed conclusion $R$, and its maps to
the original conclusion and the actual precursor contraction threads.
In particular, every source is the actual terminal premise focus of a
thread containing contraction; its frozen residuals are different
occurrences, all frozen before that focus. The thread terminates at the
focus. Each newly produced principal component starts a thread, each
unchanged background component continues its own, and no later continuation
creates a residual with that old source. Terminal references are nonempty.
In addition to these histories, the following occurrence data are retained.
\begin{enumerate}
\item The root is least: $u\le v$ means that $u$ lies on the path from $v$
to the root. Principal occurrences on different premise edges are distinguished
positions; subdivide edges with inert markers if necessary. A finite set
$\mathcal N$ of these positions is the set of residual sources.
\item The bare characteristic variable never occurs under a connective.
Each label has one left and one right occurrence at every closed node in
which it occurs. Every binary rule has label-disjoint closed premises.
Existing slots are transmitted exactly once, on the same side, through each
retained rule. A product rule groups only ordinary formulas, and a residual
rule makes only ordinary formulas principal. Unit cut cuts the constant $t$,
never a live label. Weak exchange cyclically moves the two specified word
blocks and preserves the internal order of each.
\item Every component in the root range carries a terminal-thread certificate
or a residual-source certificate. A residual of source $s$ has the ordered
slot skeleton of the designated focus at $s$. Its frozen occurrences are
carried as background after the point at which their associated live thread
reaches that focus. A retained continuation of a selected source occurrence
can meet only source positions on its path to the root. A frozen occurrence
imported from an auxiliary premise is either wholly in that imported subtree,
or has its designated future source on the rootward continuation from its
attachment position. Thus a thread cannot acquire a source above a position
at which it first enters a continuation.
\item Each terminal skeleton refers to an actual component of $G_0$.
Its slots come from the bare positions of that component, so it has at most
one slot. Ordinary formulas contain
no characteristic variable. A left slot may be filled by a live $p_a$ or
$\top$; a right slot may be filled by a live $p_a$ or $\bot$.
The root range consists of all these permitted decorations of the certified
ordered skeletons, with arbitrary finite repetition. The position-specific
skeleton at an internal proof node need not itself be in the root range.
\item Internal identities have the form $p_a\Rightarrow p_a$, with a paired
fresh label. Any slot introduced by a nonidentity initial context has its
fixed constant decoration. There is no external weakening in $\tau$.
Contraction residues have been frozen rather than contracted, so the tree
contains no external contraction. Whenever extraction imports an unselected
premise subtree, it imports its complete proof with a uniform fresh renaming.
\end{enumerate}
The component threads and source certificates are the records produced by
freezing contraction copies in preprocessing; they record a designated future
focus, not merely the position at which a formula was born.
\end{definition}

Proposition~\ref{op:preparation} supplies these data, pruning weakening
before expanding contractions and retaining slot marks through constant replacement.

\section{Separation of residual sources}\label{sec:ordered-separation}
\label{ols:section}

We prove the separation result needed for density elimination.  The proof uses
closed inputs throughout.  Its induction is on the number of selected source
positions and, for a fixed number, on a vector counting residual occurrences.
A finite support argument first removes sources shared by two or more
selected branches; a decreasing occurrence vector then removes the remaining
private sources.

\subsection{Decorated inputs and the strengthened statement}

Throughout this appendix fix an actual preparation $\tau$ as in
Definition~\ref{ols:prepared}, with at most one bare $p$ in each component
of the original end hypersequent $G_0$ before pruning. The original
conclusion and its precursor history are fixed throughout the recursion.  The tree is finite, its root is the least position, and its contraction
threads have designated terminal positions.  A thread terminating at a genuine
focus gives a \emph{source}; a thread terminating at a root component gives a
\emph{terminal component}.  If necessary, a unary marker on a premise edge
distinguishes its focus position from the other premise of the same binary
rule.  Thus the source set $\mathcal N$ is a finite set of positions in a rooted
tree, not a set of formula shapes.  We write $u\leq v$ when $u$ lies on the path
from $v$ to the root, and $u\parallel v$ when the two positions are incomparable.

At the root-range interfaces used below, every component occurrence carries its designated source or terminal-thread certificate.  A frozen residual of source $s$ has the slot skeleton $\sigma_s$.
The source records the future focus of its contraction thread; it need not be
the position where the component's formulas were first formed.  Frozen
residuals subsequently travel as background.  Every terminal skeleton has at
most one marked slot.  Ordinary unmarked formulas contain no $p$.

A slot is a specified occurrence of the original bare variable $p$, including
its side and its identity in the skeleton.  Its permitted decorations are
$p_a$ or $\top$ on the left, and $p_a$ or $\bot$ on the right.  Unmarked
occurrences of $\top$ and $\bot$ are ordinary formulas and are never changed by
slot instantiation.  A component is in the \emph{decorated range} if it is an
allowed decoration of its certified source or terminal skeleton.  There may be
arbitrarily many copies of a certified component.  In particular, belonging to
the decorated range does not mean having the same actual formula as one fixed,
already decorated component at the root of $\tau$.

The prepared rules have the following occurrence property.  Every existing
slot is carried to exactly one occurrence on the same side of the conclusion.
No slot is a logical principal formula.  Binary premises have disjoint live
labels.  An internal $p_a\Rightarrow p_a$ introduces a paired, fresh live label;
an internal slot originating in a weakening or an initial context has its
fixed, allowed constant decoration.  We use the generalized lattice rules
from preprocessing, so their two premise contexts are also transmitted once.
The tree contains neither external weakening nor external contraction; cuts, if present, are only unit cuts.  These
facts allow extraction along actual component occurrences.

\begin{definition}[Closed components and strict copies]
\label{ols:closed}
A decorated hypersequent is \emph{closed} if every live label occurs exactly
once on the left and exactly once on the right.  Its label graph has component
occurrences as vertices and live labels as edges; the endpoints of the edge
$a$ are the components containing its left and right occurrences.  Loops and
parallel edges are allowed.  A \emph{closed unit} is a connected component of
this graph, including a label-free isolated vertex.

Two closed blocks are strict copies if a single bijection of their live labels
and a bijection of component occurrences identify their actual ordered antecedents and succedents,
source or terminal certificates, slot positions and decorations, and current
protection marks.  The rule $\ECs$ deletes one of two such disjoint closed
blocks.  It never identifies two vertices or two edges.  Copy-history
identifiers are transported by the copy correspondence; they need not be
literally equal.
\end{definition}

The graph definition agrees with minimal nonempty closed subhypersequents:
a subset of vertices is closed precisely when no edge has exactly one endpoint
in it.  Constant-filled slots contribute no edges.  Consequently a terminal
vertex has degree at most one, even after its slot decoration is changed.
Strict copying compares actual formulas, so it never identifies a live $p_a$
with $\top$ or $\bot$ merely because they occupy the same slot.

\begin{definition}[Fresh-copy inputs and protection]
\label{ols:conventions}
A derivation from closed inputs $X_1,\ldots,X_m$ may use finitely many leaves,
each a uniform fresh renaming of one complete $X_i$.  Only live labels are
renamed; actual constant decorations and all certificates are copied.  This is
a convention about the permitted hypothesis leaves, not an additional
inference from a fixed labelled $X_i$ to every renaming of $X_i$.

Residual occurrences carry finite stacks of protection marks.  A nonempty
stack makes an occurrence \emph{passive}: it is never selected as an elimination
focus.  Marks already present when a procedure is called are permanent for
that call.  The call may push temporary marks and must remove exactly its own
marks before returning.  Copies inherit the stacks, and strict contraction
preserves them.  A passive occurrence is only carried as background, copied
with a whole input or proof, or deleted within a legal strict contraction.
\end{definition}

Whenever an input already has a proof without hypotheses, each of its allowed
leaves has such a proof: apply one uniform fresh bijection to the entire
original proof.  All internal identities, rule instances, closedness and
premise-disjointness conditions are preserved.  Substitution of these proofs
at the finitely many hypothesis leaves therefore produces an ordinary
hypothesis-free proof in $\Om$.

For a nonempty antichain $I\subseteq\mathcal N$, put
\[
 \downarrow I=\{s\in\mathcal N:(\exists i\in I)\ s\leq i\},
 \qquad
 \Par(I)=\{s\in\mathcal N:(\forall i\in I)\ s\parallel i\}.
\]
Relative to $I$, an unprotected residual whose source is not in $\Par(I)$ is
called \emph{bad}.  Already good residuals may be temporarily protected; their
marks are removed on return.  This convention does not remove them from the
final conclusion.

\paragraph{Key invariants for separation.}
The construction below keeps the following three occurrence invariants explicit.
\begin{enumerate}[label=(\roman*)]
\item \emph{Actual source histories.}
Each residual-source certificate refers to an actual precursor contraction
thread. Its original residuals were frozen before the designated live focus
was used. Each newly formed principal component starts a new thread, and
each unchanged background component continues its original thread. Later copies retain
this source history, as supplied by Proposition~\ref{op:preparation} and
Definition~\ref{ols:prepared}.
\item \emph{Passive decorations and copy histories.}
While a residual remains passive, its surviving descendants retain its
ordered formulas and slot decorations, up to uniform fresh renaming of a
complete input or proof. Each fresh copy has its own history record and
retains its correspondence with the original occurrence. Strict contraction
deletes whole closed copies; it does not merge their label graphs
(Definitions~\ref{ols:closed} and~\ref{ols:conventions}).
\item \emph{Terminal degree at most one.}
Every terminal skeleton in the decorated root range has at most one marked
slot. Its vertex therefore has degree at most one, under every permitted
decoration, and cannot serve as a junction between different label families.
The further conclusion that the relevant closed units lie within one family
or one private-pool copy also uses the residual-source analysis below.
\end{enumerate}
These invariants summarize the prepared occurrence data and their treatment
under copying and protection; they add no hypotheses to the following theorem.
Their preservation and their role in the family and private-pool bounds are
verified in Lemmas~\ref{ols:families} and~\ref{ols:pools} and in the proof below.

\begin{theorem}[Separation with arbitrary passive background]
\label{ols:main}
For the preparation of Definition~\ref{ols:prepared}, let $I\subseteq\mathcal N$ be a nonempty finite antichain.  For every $i\in I$
let $X_i$ be a closed hypersequent in the decorated range.  Its unprotected
residuals have sources at most $i$, and it contains at least one unprotected
residual of source $i$.  Its passive residuals may have arbitrary sources.

There is a finite derivation from these inputs, in the sense of
Definition~\ref{ols:conventions}, to a closed hypersequent $F$ in the decorated
range such that every residual of $F$ is either a copied descendant of a
residual passive at the start of this call, or has source in $\Par(I)$.
All marks introduced by the call are removed before it returns.  In particular,
no newly produced source comparable with $I$ is concealed under a permanent
protection mark.

The derivation consists of instantiated extraction templates and $\ECs$.
All boundaries of these macro steps remain in the decorated root range;
the internal nodes of an expanded template instead carry the appropriate
position-specific slot certificates and form a safe $\Om$ derivation.
Every elimination focus lies in $\downarrow I$; each joint template has at most
$|I|$ open focus leaves; and the new residuals supplied by any such template,
as opposed to its carried input background, have no source $s\geq i$ for any
$i\in I$. Passive input occurrences retain their actual decorations.
The construction also returns the occurrence and label records specified
before Lemma~\ref{ols:families} for this same derivation.
\end{theorem}

We prove this strengthened assertion for every choice of passive backgrounds.
The family and occurrence records in the following lemmas will be maintained
in the construction; they are consequences of the rule-by-rule construction,
not further assumptions on the inputs.

\subsection{Instantiation, extraction and closed background}

\begin{lemma}[Open slot instantiation]
\label{ols:slots}
Let an extracted template have designated open leaves with slot skeletons
$S_1,\ldots,S_w$.  For each leaf provide a complete closed input
$B_a\mid S'_a$, where the designated occurrence $S'_a$ has skeleton $S_a$.
After uniformly freshening the complete inputs to make their live labels
pairwise disjoint, assign to each boundary slot the actual decoration of its
corresponding slot in $S'_a$.  Give every complete auxiliary proof fresh live
labels, retaining its actual slot decorations.  The resulting template, with
each $B_a$ carried as background, is a legitimate closed derivation.
\end{lemma}
\begin{proof}
Boundary slots are occurrences, not equations between labels.  In particular,
two boundary slots that happened to have the same label in the original tree
are instantiated separately.  The actual complete input decides whether their
live decorations coincide.  A slot filled by a constant contributes no live
label.  The two occurrences of an internal identity $p_b\Rightarrow p_b$,
by contrast, receive the same fresh label.

Trace each boundary slot along the retained rules.  It never becomes a logical
principal formula, changes side, splits or disappears.  Replacing all of its
successors by its assigned decoration consequently gives another instance of
each retained rule.  This observation includes a change from a constant-filled
boundary slot to a live slot.  Only marked occurrences on this boundary trace
are changed: neither an ordinary constant nor an internal auxiliary slot is
altered. Any slot in a nonidentity auxiliary initial context keeps its
fixed extreme-constant decoration, so that the initial rule remains legal.

For a live input label, its two occurrences are either both in the transmitted
fragment, one in that fragment and one in $B_a$, or both in $B_a$.  Every existing
occurrence is carried once.  Thus its two endpoints remain present and
same-side uniqueness is preserved at every node.  An internal label has its
two endpoints in one complete fresh auxiliary proof.  At a binary combination,
the complete incoming subderivations have disjoint label sets.  These
observations prove closedness and disjointness by induction over the finite
template. There is no internal exchange hidden in this argument. For
communication, replace each slot in its specified ordered block, and exchange
exactly the two blocks selected by the original rule; their internal words
are unchanged. For weak exchange, perform the same prescribed cyclic block
rotation. For unit cut, the cut formula remains $t$ and all live slots remain
in the untouched context words. Generalized meet and join transmit each
premise context once, so they introduce no duplicate slot occurrences.  They also prove that every retained background occurrence keeps its
actual decoration and certificate.  The complete output, including its retained external background, is closed.
\end{proof}

\begin{lemma}[Lifting a derivation through a closed background]
\label{ols:context}
Suppose a finite derivation $\mathcal D$ derives $F$ from closed inputs
$U,V_1,\ldots,V_t$ and actually uses an occurrence of the hypothesis $U$.
For a closed block $C$ disjoint from the relevant labels, there is a derivation
of $C\mid F$ from $C\mid U,V_1,\ldots,V_t$.  The components of $C$ remain
background throughout. Contraction here is the coherent closed-block rule $\ECs$.
\end{lemma}
\begin{proof}
At every used $U$ leaf put a complete fresh copy of $C$ beside the corresponding
copy of $U$.  Reproduce $\mathcal D$, carrying the added components as background.
Each old inference has the same specified principal occurrences as before.
An old strict contraction acts on its old closed blocks; the newly added $C$
copies are not included in them.  Uniformly freshen proof copies before every
binary combination.

At the root there are $F$ and finitely many whole copies of $C$.  Since no rule
has acted on these copies, their actual formulas, decorations, certificates
and protection stacks correspond under a single uniform renaming of each
copy.  Contract them to one by $\ECs$.  If a proposed derivation does not use
$U$, its proof from the remaining inputs can instead be used directly whenever
its conclusion already meets the required goal.  Neither alternative projects
$U$ from $C\mid U$.
\end{proof}

Extraction is defined on the undecorated occurrence tree and then instantiated
by Lemma~\ref{ols:slots}.  Starting at a designated focus, retain that component
and follow the path to the root.  A retained component that is background at a
rule is simply carried past that rule.  If a retained principal component is
used by a unary rule, replay the rule.  If it is used by a binary rule, import
the complete other premise and replay the rule, unless a retained fragment
from that premise is already available.  Every component retained in this way
ends at its certified terminal or frozen-residual occurrence.  Thus outputs
of complete templates are in the decorated range, although their internal
nodes need not be root-component skeletons.

\begin{lemma}[Ordered source exclusion and pruned joint extraction]
\label{ols:joint}
Choose one focus occurrence for each source in an antichain $I$.  There is a
finite extraction template using each selected occurrence in a nonempty
subset $J\subseteq I$ exactly once.  Every new residual in its output has
source $s$ satisfying
\begin{equation}
             i\not\leq s\qquad\hbox{for every }i\in I.       \label{ols:exclusion}
\end{equation}
This assertion concerns the entire original $I$, including unused selections.
It is unaffected by slot instantiation.  For $I=\{i\}$ it gives the ordinary
single-focus elimination template.
\end{lemma}
\begin{proof}
Construct retained fragments from the selected positions towards the root.
For a unary rule and for a binary rule with selections in only one premise,
use the procedure just described.  At a binary rule whose two premise
subtrees both contain selections, two retained fragments have already been
constructed.  If both contain the respective principal occurrences, replay
the rule with these fragments.  If one fragment lacks its principal
occurrence, keep that fragment, discard the other fragment, and skip the rule:
every component of the kept fragment is background at that rule.  If both
lack the principal occurrence, choose either fragment.  These operations keep
a nonempty set of selected leaves.  A discarded fragment, its auxiliary
proofs and all its previously introduced occurrences are discarded together.

For a selection that survives, its open leaf omits the other components
at the selected source position, including all frozen copies of that
thread. Every residual with source $s\geq i$ was frozen before its actual
source focus; hence it is already frozen on or above the cut at $i$.
It cannot be newly created below $i$: the thread ending at $s$ stops at
that focus. Each newly formed principal component there starts a new
thread, while each unchanged background component retains its own thread. In particular,
the continuation cannot recreate a residual of the selected source $i$.  Nor can
it enter through an auxiliary premise: a thread whose terminal focus is on
or above $i$ lies in the already cut subtree, whereas a later auxiliary
premise is a disjoint subtree.  This argument concerns component threads and
their designated terminals, not the position at which a formula first arose.

If a selection $i$ is discarded at a binary node, every position on or above
$i$ is inside the discarded premise subtree.  Below that node, this entire
subtree remains inside the already traversed premise.  A subsequently
imported other premise cannot meet it.  Therefore no source on or above $i$
can re-enter later.  Applying these two arguments to every original selection
proves the exclusion statement for all of $I$.

Every replayed inference is an original rule restricted to its retained
principal occurrences and carried background.  Skipping a rule keeps only
components that it transmits unchanged.  The remaining non-open leaves are
complete auxiliary proofs.  Lemma~\ref{ols:slots} therefore turns the result
into a legitimate derivation from complete closed inputs.  Decoration changes
neither these source positions nor any pruning decision.
\end{proof}

\paragraph{Occurrence, family and copy records.}
The following records distinguish the origins of labels from the groups used
later to count private occurrences. They are produced by the construction,
not additional conditions on an input. An input may already carry records
from an earlier call; a new call retains those records beneath its own
whole-input copy record.

For a call with inputs $(X_i)_{i\in I}$, distinguish its used whole-input
copies by identifiers $d$, write $t(d)$ for the input index, and write
$\rho_d$ for the uniform renaming of that copy. To compare ranges, use
disjoint tagged copies of the label alphabet, embedded in a common fresh
alphabet. On the finite support of any actual copy its embedding is a
bijection onto the labels of that copy. For each output occurrence retained under
an input protection mark, record $(d,b)$, where $b$ is its position in
$X_{t(d)}$. This map is injective, and the recorded occurrence is the
actual component at $b$ under $\rho_d$, including its slots and its
inherited marks after the call's temporary marks are removed.

Each live label has an origin record of one of two kinds. An input origin
$(d,\ell)$ records the label $\ell$ in $X_{t(d)}$ and has actual value
$\rho_d(\ell)$; its \emph{input family} is $d$. A label introduced by a
complete fresh auxiliary proof has an \emph{auxiliary family} $a$, the
identifier of that proof import, together with its own label identity.
These identifier ranges are disjoint. Conversely, a live output label in
the range of $\rho_d$ must have the input origin $(d,\ell)$ for an actual
live label $\ell$ of $X_{t(d)}$. The label record agrees with the position
record on each retained component. A call also records a common auxiliary
family for all live slots of each residual that it returns active; the
construction of this part of the record is proved in the next lemma and
in the simultaneous induction for Theorem~\ref{ols:main}.

Composition retains the finer origin record. If a fresh copy $e$ of an
intermediate output contains a position previously recorded as $(d,b)$,
its composite record is $((e,d),b)$, and the corresponding input-label
origin is $((e,d),\ell)$. An inherited auxiliary family $a$ becomes
$(e,a)$; a new auxiliary import receives a separate fresh identifier.
Thus regarding an entire intermediate output as one input of a recursive
call does not identify its earlier families. The local whole-input
identifier and the finer ancestry record are both retained. Strict
contraction deletes records of the discarded closed block and restricts
all maps to the surviving positions and labels. No contraction identifies
two surviving label origins. These records alone assert neither a
private-unit bound nor a decrease of the occurrence measure.

\begin{lemma}[Fresh families in the separation construction]
\label{ols:families}
At a joint extraction, or a root-range stage of the construction of
Theorem~\ref{ols:main} obtained by composing extractions and recorded
recursive returns, every newly supplied active residual with source in
$\Par(I)$ has a common auxiliary identifier for its live labels, relative
to the inputs of the call. Here every focus lies in $\downarrow I$ and
the only active input sources lie in $\downarrow I$. A surviving copy of
such a residual keeps a common auxiliary identifier when its temporary
protection is removed.

At a completed return, every closed unit containing descendants of
initially passive occurrences contains such descendants from at most one
whole-input copy. The assertion includes initially passive occurrences
with source in $\Par(I)$ and label-free isolated occurrences.
\end{lemma}
\begin{proof}
In an extraction assign fresh, disjoint names to the open input copies
and the imported auxiliary proofs. Each surviving passive position has
one input position. An occurrence newly formed on the continuation of
a selected focus $u$ can have a source only on its path towards the root,
hence at most $u$. If the thread of a frozen component from an auxiliary
premise has its source outside that premise, its source likewise lies
on the receiving continuation (Proposition~\ref{op:preparation}).
Thus a newly supplied residual of source in $\Par(I)$ comes wholly from
one auxiliary premise and all its live labels have a common auxiliary
identifier. Transmitting slots, dropping a discarded premise, and
deleting a strict copy preserve the origins of all retained labels.
The disjoint renamings give the stated reflection of input-label origins.

For composition, substitute the records of the actual returned proofs.
A position that passes through a fresh copy $e$ of an earlier output and
then traces to $(d,b)$ has record $((e,d),b)$. The composite is injective
because both recorded maps are injective. Transport live-label origins
with the same pair of identifiers. A label recorded as auxiliary in the
earlier output receives its earlier identifier tagged by $e$; auxiliary
identifiers of the outer operation have a separate range. Reflection at
each stage ensures that an inherited input label has exactly its recorded
origin after composition. A component copied under protection keeps its
former labels. A newly active residual returned by a recursive call has
the common auxiliary identifier provided by that smaller call. An earlier
good residual copied under a temporary mark retains its former common
identifier when the mark is popped. A strict contraction restricts the
records to the surviving positions. These operations prove the first
assertion at each stage of the construction. In recursive uses the
records are supplied by the induction in Theorem~\ref{ols:main}.

At a completed return, a residual either descends from a passive input
and has the family of its recorded copy, or is active and has a common
auxiliary identifier. A terminal vertex has degree at most one by
Definition~\ref{ols:prepared}. No vertex can connect edges of different
families, so along a path every edge has the same family. The input
family of any labelled passive descendant is its actual copy identifier.
A label-free descendant forms a singleton unit, and its position record
gives the copy even in this case.
\end{proof}

\subsection{Private intervals and finite occurrence counts}

Assume $|I|>1$.  The lowest common ancestor here means the greatest common
lower bound in the root-minimal order.  Define
\begin{equation}
 h_i=\max\{\lca(i,j):j\in I\setminus\{i\}\},\qquad
 P_i=\{s\in\mathcal N:h_i<s\leq i\}.                  \label{ols:private}
\end{equation}
All the nodes in the displayed maximum lie on one root path, so the maximum
exists.  A source in $P_i$ is \emph{private to $i$}.  A source $q\in\downarrow I$
that is not private to any $i$ lies below at least two members of $I$ and is
called a \emph{bridge source}.  Indeed, if $q\leq i$ and $q\not>h_i$, then
$q\leq h_i=\lca(i,j)\leq j$ for some $j\ne i$.  Conversely a source below two
members cannot lie in a private interval.  These definitions use all of $I$,
not merely the selections retained after pruning.

Apply Lemma~\ref{ols:joint} to complete inputs $X_i$, selecting one unprotected
occurrence of source $i$ in each, and denote the actual output by $H_0$.
For every used $X_i$, make a \emph{private pool record} consisting of its whole
input copy and the complete auxiliary proofs imported along its private
path above $h_i$. The record includes component-occurrence histories as well
as its set of live labels, and is retained even if that set is empty.
Auxiliary proofs imported elsewhere receive separate records. The label
sets of distinct records are disjoint.

\paragraph{Pool colours and their transport.}
A pool is a grouping of the origin families just defined, not a new label
origin. At $H_0$, assign the input family of a used $X_i$ and every auxiliary
family imported on its private path to the same pool $i$. Give each auxiliary
import outside those paths a separate pool identifier. Every finer family inherited inside that whole $X_i$ copy receives the
same initial pool $i$. This defines a map
$\kappa_0$ from the origin families of $H_0$ to pool identifiers; every
live label receives the colour of its family under $\kappa_0$. The position
record assigns an input occurrence to its input pool even when it has no
live label. No label is renamed or identified by this grouping.

For bookkeeping that distinguishes only the private pools, one may compose
$\kappa_0$ with the map that fixes every private index $i$ and sends all
non-private pool identifiers to one neutral colour distinct from every
private index. This coarsens the
colouring without identifying origin-family records, component occurrences
or whole-block copies. Distinct private pools retain distinct colours;
the private-pool bounds and copy-position injections below use these
private indices together with the unchanged copy records.

The following distinctions will be used in the bridge execution.
\begin{center}
\begin{tabularx}{\linewidth}{@{}lX@{}}
\toprule
Record & Meaning and behaviour under copying\\
\midrule
Input family $d$ & One complete copy of an input to the current call;
its finer records from earlier calls are retained.\\
Pool $i$ in $H_0$ & The used $X_i$ copy and its private-path auxiliary
imports, grouped by $\kappa_0$.\\
Whole-$H_0$ history $d_0$ & The ancestry of one complete copy of $H_0$;
a further copy $e$ extends this history to $(e,d_0)$.\\
Transported colour $(d_0,i)$ & Pool $i$ inside that particular
whole-$H_0$ copy, rather than all copies of pool $i$ together.\\
Supplemental input family & A fresh complete original input used in a
bridge call, with a colour disjoint from the copied current state.\\
New auxiliary family & A complete fresh auxiliary-proof import, with
its own colour disjoint from input families and transported colours.\\
\bottomrule
\end{tabularx}
\end{center}
A label descended from $H_0$ retains its $H_0$ pool; other inherited
families retain their separate colours. A whole-state copy tags every
inherited colour with its copy history. Supplemental inputs and new
auxiliary families have separate identifier ranges. When recursive returns are composed, the finer origin
records determine which of these cases applies; the recursive call's
coarser input-family name never replaces an old pool colour. The proof of
Lemma~\ref{ols:private-copy-injection} establishes the resulting coloured-graph
invariant at each bridge-step boundary after the temporary marks have
been restored.

\begin{lemma}[Private-pool bound]
\label{ols:pools}
In $H_0$, every unprotected bad residual private to $i$ belongs to the private
pool of an actually used $X_i$, and all its live labels are in that pool.
Every originally passive residual belongs to the pool of its actual input,
regardless of its source.  A bad residual mixing private pools must have a
bridge source.  If $n_i$ is the number of unprotected source-$i$ occurrences
in $X_i$, the private pool contains at most $n_i-1$ unprotected occurrences
of source $i$, and no bad private occurrence above $i$.
\end{lemma}
\begin{proof}
Above $h_i$, no path from a different selected position has joined the path
from $i$.  Every rule retained there uses either the retained $i$ fragment
or a complete auxiliary proof imported on that path.  At a source position
$t\in P_i$, a residual assigned to $t$ is already frozen and thereafter
travels unchanged as background.  Trace that concrete component backwards
from $t$.  It can meet only the $i$ input or an auxiliary subtree imported
between $i$ and $t$.  It cannot contain labels from another selected branch,
or from a subtree first introduced below $h_i$.

This also covers a residual whose formulas were formed strictly before
its source position. Suppose that a residual of private source $t$ is
imported from an auxiliary subtree attached to the selected path at $v$.
An auxiliary subtree is disjoint from that selected path. Its source $t$,
which lies on the path between $h_i$ and $i$, therefore cannot be a past
source internal to the auxiliary subtree. The import dichotomy of
Proposition~\ref{op:preparation} puts $t$ on the future continuation from
$v$, so $t\leq v\leq i$. Since $t>h_i$, the attachment point $v$ is
itself in the private interval, and its complete auxiliary family is part
of the $i$ pool. A residual with a past source inside the auxiliary subtree
may also be imported, but that source is not private to $i$. This
explicit distinction is why the argument does not identify a formula's
birth with its residual source.  Private residuals carried unchanged
from the original input already have that input's labels.  Passive input
residuals also remain in their own input pool; their sources need not be
private and need not equal $i$.

Two private paths can first interact only at their common lower ancestor or
below it.  Any bad residual whose labels mix their pools consequently has a
source below both selected positions, and is a bridge residual.  If a
selection is pruned, its whole retained fragment and its introduced records
are deleted.  The discarded original subtree cannot subsequently be imported
as an auxiliary subtree, by the disjointness argument in
Lemma~\ref{ols:joint}.  No later private residual can therefore be attributed
to an unused input.

The selected source-$i$ occurrence in each used input has been replaced.
Lemma~\ref{ols:joint} supplies no new source $i$ or source above $i$.
Other input occurrences are carried only once in this template.  Hence at
most $n_i-1$ source-$i$ occurrences remain in its pool.  This count concerns
occurrences, including those whose slots are all constant-filled.
\end{proof}

For a source $s$, let $r(s)$ be the length of its rooted premise-edge
address, including the final edge identifying the premise focus. Thus
$r(s)\geq 1$ and $u<v$ implies $r(u)<r(v)$. Incomparable sources may
have the same rank. Put $N=\max_{s\in\mathcal N}r(s)$, with $N=0$ when
$\mathcal N$ is empty. For a protected hypersequent $X$, let
\begin{equation}
 \mu(X)=(c_N(X),\ldots,c_1(X)),\qquad
 c_a(X)=\#\{\hbox{unprotected residual occurrences of rank }a\}.       \label{ols:measure}
\end{equation}
Whenever a call begins, its good residuals are protected, so only its bad
sources are counted.  The strict lexicographic order on these vectors has no infinite decreasing
chain (when $N=0$ there is only the empty vector). For $N>0$ it is represented by the ordinal
$\omega^{N-1}c_N(X)+\cdots+\omega c_2(X)+c_1(X)$.  The vector counts source occurrences, not slots or labels.
Turning all slots of a residual into constants therefore does not remove its
contribution.

For a bridge source $q$ occurring unprotected in a current output $H$, put
\begin{equation}
 I_q=\{q\}\cup\{i\in I:q\parallel i\}.                 \label{ols:bridgefrontier}
\end{equation}
A bridge step uses the whole $H$ as the $q$ input. Its unprotected
residuals with source at most $q$ remain unprotected, and all its other
unprotected residuals receive a temporary mark. The remaining inputs are
fresh complete copies of the original $X_i$ for $i\in I_q\setminus\{q\}$.
After the call, its temporary marks on the surviving copies of $H$ are
removed, and the prior protection stacks are restored.

Since $q$ lies below at least two members of $I$, $I_q$ is an antichain
of smaller size. Every $i\in I$ either lies above $q$ or belongs to
$I_q$. If $q\leq i$ and $u$ is comparable with $i$, then either
$q\leq i\leq u$, or $u,q$ lie on the root-to-$i$ path and are
comparable. Hence
\begin{equation}
                          \Par(I_q)\subseteq\Par(I).          \label{ols:coverage}
\end{equation}

\begin{lemma}[Occurrence injection for a single private-pool copy]
\label{ols:private-copy-injection}
Fix the joint output $H_0$ of Lemma~\ref{ols:pools}, with its actual
position and pool records. Let $H_1$ be the endpoint of a finite bridge
execution from $H_0$ with no unprotected bridge source remaining. An
execution may be empty. Each nonempty step chooses a current unprotected
bridge source $q$, uses the actual recursive return for the frontier
$I_q$ with the occurrence and label records described before
Lemma~\ref{ols:families}, and takes its output after restoring the prior
marks. Good residuals are protected when counting bad occurrences.

For a used input $i$, let $B_i$ be the positions of unprotected residuals
of $H_0$ with source in $P_i$. For a fixed copy history $d$ of the
\emph{whole} $H_0$, let $B_{d,i}$ be the positions of surviving unprotected
bad occurrences in $H_1$ with that history and source in $P_i$. There is
a source-preserving injection
\[
                     \iota_{d,i}:B_{d,i}\hookrightarrow B_i.
\]
Every closed unit $U$ of $H_1$ containing an unprotected bad occurrence
admits a single pair $(d,i)$ for all its bad occurrences, and every live
label of $U$ has the colour of that pool copy. For this same $i$, with
$n_i$ as in Lemma~\ref{ols:pools},
\begin{equation}
 c_a(U)=c_a(X_i)=0\quad (r(i)<a\leq N),\qquad
 c_{r(i)}(U)\leq n_i-1<n_i=c_{r(i)}(X_i).
 \label{ols:private-copy-ranks}
\end{equation}
In particular,
\begin{equation}
                         \mu(U)<_{\rm lex}\mu(X_i).           \label{ols:unitdecrease}
\end{equation}
\end{lemma}

\begin{proof}
A private source cannot be at most a bridge source, since it would
then lie below two members of $I$. Its occurrences in the current $H$
therefore receive or retain protection in each bridge call. By the
classification in the recorded return and~\eqref{ols:coverage}, a bad
occurrence unprotected after the marks are restored descends from a
temporarily protected occurrence of $H$. Passive occurrences in the
supplemental inputs retain their original marks. Iterating shows that a
surviving private bad position traces to a private bad position of $H_0$.

At the empty execution use the identity on the private bad positions of
$H_0$. At a further step the retained-position injection for the child
assigns a bad position to a pair $(e,k)$, where $e$ is a whole copy of
the preceding $H$ and $k$ is a position in it. If $k$ previously traced
to $(d,b)$, assign $((e,d),b)$ to the new position. Equality of these
records forces equality first in the preceding injection and then in
the child's injection. The component equality in the records preserves
the source, including when every slot is constant-filled. A strict
contraction just restricts the resulting injection. Fixing a whole-copy
history and $i$ gives $\iota_{d,i}$ into the positions of $H_0$ itself.

To determine which injection applies to a closed unit, use the pool
colours and whole-copy histories defined above. The coloured-graph invariant
at $H_0$, and at each subsequent bridge-step boundary after restoring its
temporary marks, is:
\begin{quote}
Every residual vertex has all its incident live labels of one colour,
except possibly an unprotected bad residual with a bridge source.
A terminal vertex has degree at most one. An inherited private bad
vertex has the transported colour $(d,i)$ specified by its position
record, where $d$ is its whole-$H_0$ copy history.
\end{quote}
A label-free vertex has no incident edge; its position record is still
retained. Thus the invariant does not discard constant-filled occurrences.

At $H_0$, a private bad or initially passive vertex has one pool colour
by Lemma~\ref{ols:pools}. A newly supplied good residual has one auxiliary
family by Lemma~\ref{ols:families}, and hence one pool colour under
$\kappa_0$. Any remaining unprotected bad source is a bridge source.
This proves the initial invariant.

For a bridge step, transport an old colour with the identifier of the
whole copy and use disjoint colours for supplemental input copies and
new auxiliary families. The origin-reflection property of the recursive
return identifies old labels with their finer records, so their colours
agree with the position injection. Protected old occurrences keep their
actual decorations and hence their incident colours. In particular, old
good residuals stay protected, since their sources cannot be at most the
bridge source $q$, which lies below a member of $I$. A surviving passive
supplemental occurrence has the single colour of its supplemental input
copy. A new active residual in the return has a single auxiliary family;
by~\eqref{ols:coverage} it is good for the original $I$.

An old mixed-colour bridge residual may temporarily be protected inside
this recursive call. At the recursive root-range interfaces this protected old vertex need not
be monochromatic. When the call's marks are restored, any such surviving occurrence
again has its prior unprotected bridge status. An occurrence passive
before this bridge step's temporary marks were added remains monochromatic;
no temporary mark can make a mixed bridge residual permanently passive. These facts prove the invariant at
the next bridge-step boundary. They also show why the invariant is asserted
after mark restoration, not at arbitrary internal nodes of a recursive call.

At $H_1$ no unprotected bridge source remains, so the exception has been
eliminated. Every residual vertex is now monochromatic, and a terminal
vertex cannot connect two edges. Consecutive edges on a graph path
therefore have the same colour. A closed unit containing a labelled bad
position consequently has the colour $(d,i)$ of that position on every
edge and at every bad position. A label-free bad position forms an isolated
unit, and its position record gives the same single-pool-copy conclusion.

All counted sources in this $U$ are in $P_i$, and every counted source
of $X_i$ is at most $i$. In the chain below $i$, the only source of
rank $r(i)$ is $i$, and no source has greater rank. Restrict
$\iota_{d,i}$ to source-$i$ positions of $U$. In $B_i$ there are at
most $n_i-1$ such positions by Lemma~\ref{ols:pools}; the original
$X_i$ had exactly $n_i$. This proves~\eqref{ols:private-copy-ranks}
and gives~\eqref{ols:unitdecrease} at coordinate $r(i)$.
\end{proof}

\subsection{Proof of the strengthened separation theorem}

\begin{proof}[Proof of Theorem~\ref{ols:main}]
We use induction first on $m=|I|$, and, for fixed $m$, on the lexicographic
vector $\sum_{i\in I}\mu(X_i)$.  The induction covers all antichains of that
size and all passive backgrounds.  The induction constructs the separation derivation, the stated source
restrictions, and its occurrence and label records simultaneously. At a
recursive step we use the records already returned by the strictly
smaller call. In particular the private-unit bound will be derived from
the actual smaller-arity bridge returns before any same-arity call is made.

\paragraph{One selected source.}
Let $I=\{i\}$.  If a bad unprotected source is present, choose a greatest
such source $s\leq i$ and one occurrence of it.  Apply its single-focus
template to the whole current closed hypersequent, retaining the other
components as background.  Lemma~\ref{ols:joint} excludes every new source
$\geq s$.  A newly supplied source comparable with $i$ must be below $s$:
if it lies above $i$ it is above $s$, and if it lies below $i$ it is comparable
with $s$ on the root-to-$i$ chain.  All new incomparable residuals are good
and are temporarily protected.

The chosen source-$s$ occurrence disappears. All counted sources lie on
the chain below $i$, so equal rank there implies equal source. Higher-rank
counts are unchanged and every newly counted occurrence has lower rank.
The vector thus strictly decreases.  The process terminates, preserves all
previous protection marks, and leaves only their descendants and sources
in $\Par(i)$.  Remove the marks introduced for good residuals.  Every focus
was at most $i$, and each macro excluded sources above $i$ because its focus
was at most $i$. Lemma~\ref{ols:slots} preserves the range. At each
application record the surviving passive positions, take disjoint names
for the imported auxiliary proofs, and compose with earlier copy records
as in Lemma~\ref{ols:families}. Popping a good residual's mark preserves
its recorded family.

\paragraph{Initial joint step.}
Suppose $m>1$ and all smaller-arity assertions are known, independently of
their occurrence counts.  For the present arity assume the assertion for
strictly smaller total vectors.  Apply the pruned joint template to obtain
$H_0$ as above.  There is at least one used input.  Its new sources exclude
all positions on or above any member of the original $I$.  Therefore every
unprotected bad source in $H_0$ lies in $\downarrow I$.  Temporarily protect
newly produced good residuals. Retain the position and live-label
origins of this extraction, together with the private-pool assignment
of Lemma~\ref{ols:pools}.

\paragraph{Removal of bridge sources.}
Whenever the current output $H$ has an unprotected bridge source $q$,
use the inputs for the frontier $I_q$ in~\eqref{ols:bridgefrontier}.
The $q$ input has an unprotected occurrence of source $q$ and all its
unprotected sources are at most $q$. The other inputs have the required
source bounds by hypothesis. Since $|I_q|<|I|$, the outer induction
supplies this very separation derivation and its occurrence and label
records. Append its recorded return to the bridge execution, composing
its copy identifiers with those of the preceding state.

Remove only the temporary marks added to $H$ for this call.  Any remaining
bad unprotected occurrence is now a descendant of an occurrence of $H$
that was temporarily protected.  Bad active occurrences of the auxiliary
control inputs $X_i$ have all been eliminated, and new incomparable
residuals are good for the original $I$ by~\eqref{ols:coverage}.
Consequently no new bad source type has appeared, and type $q$ has
disappeared from the unprotected support.  Permanently passive occurrences
of source $q$ are irrelevant to this support and are allowed to remain.
Repeat while an unprotected bridge source is present.  The finite set of
bad source types strictly loses an element at each repetition, so this
procedure terminates, even though a call can copy arbitrarily many
components.  Denote its output by $H_1$.

Including the empty case, these recorded returns form the bridge
execution required by Lemma~\ref{ols:private-copy-injection}. For each closed unit $U$ of
$H_1$ containing an unprotected bad occurrence, that lemma supplies a
used input $X_i$ and the strict comparison~\eqref{ols:unitdecrease}.

\paragraph{Elimination of the remaining private units.}
Decompose the actual closed $H_1$ into its finitely many closed units.
Choose a unit $U$ containing an unprotected bad source, in an $i$ pool,
and let $k$ be its greatest such source.  The private sources form a
chain, so $k$ exists.  It satisfies $k\leq i$ and is incomparable with
every $j\in I\setminus\{i\}$.  For if $k\leq j$ it would be a bridge
source, while $j\leq k\leq i$ contradicts the original antichain property.
Thus
\[
                       K=\{k\}\cup(I\setminus\{i\})
\]
is an antichain of size $m$.  Use $U$ as its $k$ input, retain its existing
protection marks, and temporarily protect its good residuals.  Use the
original $X_j$ as the remaining inputs. Replacing the summand
$\mu(X_i)$ by $\mu(U)$ strictly decreases the total occurrence vector
by~\eqref{ols:unitdecrease}, since all other summands are unchanged.
The inner induction
therefore supplies a finite actual derivation, its retained-position
and live-label records, and the conclusion that every unprotected final
source belongs to $\Par(K)$.  Since $k\leq i$, the same root-path
argument as in \eqref{ols:coverage} gives $\Par(K)\subseteq\Par(I)$.
After removing this call's temporary marks, every residual it retains is
therefore either an originally passive descendant or good for $I$.

If this derivation ignores its $U$ input, it already proves the required
goal from the other original inputs, and the present construction may
finish with that proof.  Otherwise write the entire current output as
$C\mid U$, where $C$ is the union of all other current closed units.
Apply Lemma~\ref{ols:context} to lift the inner derivation from the actual
input $C\mid U$.  Substitute the already constructed proof of that
whole input at its leaves.  The complete background copies contract
legally to one $C$, and the selected bad unit has been replaced by an
output with no new bad unprotected residuals. In the unused-input case
keep only the records of the used inputs. In the lifted case extend the
records by the whole-copy identifiers of the substituted proofs, preserve
the records on the copied background, and restrict them to survivors
at each strict contraction. This is the composition in
Lemma~\ref{ols:families}.

Every other current unit is untouched by this lifted operation.  Perform
it for each of the originally remaining bad units.  There are finitely
many such units, and no step adds a bad unit to the unprocessed background.
This completes their elimination.  This is replacement inside an actual
closed derivation, never an inference extracting $U$ from $C\mid U$.

\paragraph{Induction invariants and termination.}
The initial joint step satisfies the asserted exclusion for the whole
original $I$.  In either kind of recursive call, every new frontier
position lies below an original member and the frontier covers every
original member from below.  A template in that call cannot newly supply
$s\geq i$ for an original $i$: choose its frontier member $a\leq i$;
then $a\leq s$, contrary to the recursive exclusion property.  Its foci
also remain in $\downarrow I$.  Bridge calls have smaller arity, private
unit calls have the same arity, and pruning only reduces the number of
open leaves.  Thus no template has more than $m$ focus leaves.

All copies and instantiated rules preserve the decorated range at the
macro boundaries by Lemma~\ref{ols:slots} and occurrence tracing.
No permanent input mark was removed, and every temporary mark has been
removed on return.  New good occurrences may have been temporarily
protected, but a new bad occurrence cannot escape either the bridge
procedure or the private-unit induction under such a mark.  The final
classification is exactly the one in the theorem. Composing the records
from the joint extraction, bridge execution and lifted private-unit
returns gives the asserted records for this completed derivation.
The connected-unit assertion of Lemma~\ref{ols:families} now applies.

Every inner recursive call has a strictly smaller finite rank vector;
every outer recursive call has smaller arity.  Between recursive calls,
the bridge loop decreases a finite source-type support, and the final
unit loop has a finite list and does not replenish it.  All extraction
templates and proof substitutions are finite.  These are well-founded
termination arguments for the whole construction, not bounds on proof
size.  The simultaneous induction is complete.
\end{proof}

\begin{corollary}[Branch separation]
\label{ols:corordinary}
Let $I\subseteq\mathcal N$ be a nonempty antichain, and for every $i\in I$
let $B_i$ be a closed hypersequent in the decorated range with
\[
 \supp_c(B_i)\subseteq\downarrow\{i\}\cup\Par(I),
 \qquad i\in\supp_c(B_i),
\]
where $\supp_c$ denotes the set of residual source positions occurring in
the hypersequent.  From uniform fresh copies of the complete inputs there
is a finite derivation of a closed $F$ in the decorated range with
$\supp_c(F)\subseteq\Par(I)$.  If the $B_i$ are provable in $\Om$, so is $F$.
The focus and source-exclusion restrictions of Theorem~\ref{ols:main} hold.
\end{corollary}
\begin{proof}
Protect exactly the input residuals whose sources already belong to
$\Par(I)$ and apply Theorem~\ref{ols:main}.  Remove these initial marks at
the end.  Their descendants remain good, and every other final residual
is good by the theorem.  When the inputs are provable, replace every used
fresh-copy leaf by the correspondingly renamed complete proof, as
explained after Definition~\ref{ols:conventions}.
\end{proof}

\begin{remark}[Normalization and the scope of the theorem]
\label{ols:normalization}
One may finish by repeatedly contracting any two strict closed copies.
Each contraction strictly decreases the number of component occurrences,
so a fully contracted conclusion exists and has the same source-range
property.  Neither uniqueness of this conclusion nor full contraction
at every intermediate stage is needed.  The theorem asserts a finite
separation derivation and its source restrictions; it does not assert
the skeleton or scheduling conditions of the former grafting argument.
\end{remark}

\subsection{Removing all residual sources}
\label{olg:section}

We now pass from the branch-separation theorem to a source-free theorem, and then recover the unnumbered strong-density conclusion.
The argument uses neither a prescribed intermediate contraction schedule nor a special skeleton labelling.
All invocations of separation are applied to entire closed, provable hypersequents.

Fix the finite numbered proof supplied by Definition~\ref{ols:prepared}, and let $\mathcal N$ be its finite set of residual-source nodes.
The proof root is least in the tree order.
For a hypersequent $H$ in its allowed slot range, define
\[
 \supp_c(H)=\{s\in\mathcal N:\text{$H$ has a residual occurrence with source $s$}\}.
\]
Terminal components have no residual source and do not contribute to this set.
For $I\subseteq\mathcal N$, put
\[
 \Par(I)=\{u\in\mathcal N:\ u\parallel i\text{ for every }i\in I\},
 \qquad M=\operatorname{Max}(\mathcal N).
\]
The allowed range includes all permitted fillings of the terminal and source slot skeletons; it is not restricted to the particular constant fillings at the preprocessed root.

\subsection{An induction over globally maximal sources}

\begin{lemma}[Provable weak branches]
\label{olg:weak-branches}
Let $I\subseteq M$ be nonempty. For each $i\in I$, suppose that $B_i$ is closed, provable in $\Om$, lies in the allowed range, and satisfies
\begin{align}
 \supp_c(B_i)&\subseteq\downarrow i\ \cup\ \Par(I),
 \label{olg:weak-condition-one}\\
 \supp_c(B_i)&\subseteq\Par(I\setminus\{i\}).
 \label{olg:weak-condition-two}
\end{align}
The occurrence of source $i$ itself is not required.
Then a closed, provable $Q$ in the same allowed range exists with
\[
 \supp_c(Q)\subseteq\Par(I).
\]
\end{lemma}
\begin{proof}
If some $B_i$ contains no source in $\downarrow i$, then \eqref{olg:weak-condition-one} already gives the conclusion with $Q=B_i$.
Otherwise, for each $i$ choose the greatest source $k_i\le i$ actually occurring in $B_i$.
This greatest source exists because all nodes below $i$ lie on one finite root path.
Mark the occurrences with sources in $\Par(I)$ passive at the input interface.
The remaining sources of $B_i$ are all at most $k_i$, and at least one occurrence has source $k_i$.

The set $K=\{k_i:i\in I\}$ is an antichain.
Indeed, \eqref{olg:weak-condition-two} implies $k_i\parallel j$ for $i\ne j$.
If $k_i\le k_j\le j$, this gives a contradiction; the reverse comparison is excluded symmetrically.
Apply Theorem~\ref{ols:main} to these actual branches and the frontier $K$.
Its final residuals are either descendants of the passive input occurrences or have sources in $\Par(K)$.
The former have sources in $\Par(I)$.
For the latter, observe that
\begin{equation}
 \Par(K)\subseteq\Par(I).
 \label{olg:frontier-cover}
\end{equation}
To prove this, let $k_i\le i$ and suppose $u$ is comparable with $i$.
If $i\le u$, then $k_i\le u$.
If $u\le i$, then $u$ and $k_i$ lie on the same root path and are comparable.
Thus $u\parallel k_i$ implies $u\parallel i$.

The separation derivation uses finitely many whole-input uniform fresh copies.
Since every $B_i$ is already a theorem, uniformly rename its entire proof for each required copy and insert that proof at the corresponding leaf.
This yields a genuine proof with no assumptions.
No inference from a fixed-labelled assumption to an independently rewired copy is involved.
\end{proof}

\begin{theorem}[Source elimination]
\label{olg:source-free}
Let $\tau$ be obtained by preprocessing an actual unit-cut proof as in
Definition~\ref{ols:prepared}. Assume that each component of its original
end hypersequent $G_0$ has at most one bare $p$ before any pruning.
For every subset $I$ of the globally maximal source set
$M=\operatorname{Max}(\mathcal N)$, one can construct a closed, nonempty
$\Om$-provable hypersequent $Q_I$ in the allowed ordered root range with
$\supp_c(Q_I)\subseteq\Par(I)$. In particular, from the prepared proof
one can construct a closed, nonempty $\Om$-provable hypersequent $H$ with no
residual-source occurrences.
Every component of $H$ is a permitted filling of a terminal slot skeleton, with arbitrary finite repetitions allowed.
\end{theorem}
\begin{proof}
For each $I\subseteq M$, construct a closed theorem $Q_I$ in the allowed range such that
\begin{equation}
 \supp_c(Q_I)\subseteq\Par(I).
 \label{olg:q-invariant}
\end{equation}
Induct on $|I|$. For $I=\varnothing$, use the root of the preprocessed proof.

For nonempty $I$, take $B_i=Q_{I\setminus\{i\}}$ for each $i\in I$.
The induction hypothesis gives \eqref{olg:weak-condition-two}.
Because $i$ is maximal in the entire source set $\mathcal N$, no source is strictly above $i$.
A source occurring in $B_i$ that is comparable with $i$ must therefore be at most $i$.
If it is incomparable with $i$, \eqref{olg:weak-condition-two} places it in $\Par(I)$.
This proves \eqref{olg:weak-condition-one}, so Lemma~\ref{olg:weak-branches} supplies $Q_I$.

Take $I=M$. Every node of a finite nonempty partially ordered set is at most a maximal node.
Hence no source in $\mathcal N$ belongs to $\Par(M)$, and \eqref{olg:q-invariant} shows that $Q_M$ has no residual source.
Set $H=Q_M$. If $\mathcal N=\varnothing$, the preprocessed root already has this property.
The construction consists of actual finite proofs, whose nodes are nonempty hypersequents; in particular it does not manufacture an empty provable branch.
Allowed-range preservation now leaves only terminal slot skeletons.
\end{proof}

The use of globally maximal sources is essential to this induction.
It ensures that each actual replacement frontier satisfies $k_i\le i$, which is the direction needed in \eqref{olg:frontier-cover}.
No replacement by a source strictly above $i$ is used.

\begin{remark}[Relation to the original separation lemma]
Theorem~\ref{olg:source-free} supplies the source-exclusion conclusion needed
from Lemma~7 of \cite{Wang2019}. It applies to every subset of the maximal
source set, including all complete subsets in the original construction.
Its well-founded recursion uses source antichains and occurrence counts;
no uniqueness of a separated fragment or prescribed intermediate
contraction schedule is needed.
\end{remark}

\section{An explicit semilinear Hilbert presentation}\label{sec:hilbert-explicit}

This appendix specifies the Hilbert target instead of leaving the
Hilbert--hypersequent correspondence implicit. Write $A\to B=A\backslash B$
and $A\rightsquigarrow B=B/A$, and put
\[
 \lambda_C(A)=(C\to A\cdot C)\wedge t,
 \qquad \rho_D(A)=(D\rightsquigarrow D\cdot A)\wedge t.
\]
The following axiom schemata define $\mathrm H_{\mathrm{sl}}^{*}$
together with the rules below. All metavariables range independently
over formulas.\footnote{For the base calculus $\mathrm{HpsUL}$, see
Metcalfe, Olivetti and Gabbay~\cite[Section~9.3, p.~243]{MOG2009}.
Definition~1 of~\cite{Wang2018} displays the same subscript $\chi$ in
both conjugates, whereas its Definition~4(iv) quantifies their parameters
independently. Here we adopt the independent-parameter specification,
consistent with Definition~4(iv), and denote the resulting calculus by
$\mathrm H_{\mathrm{sl}}^{*}$. We do not assume or claim equivalence
with the calculus obtained by restricting PRL to a single shared
parameter.}
\begin{align*}
&\text{(A1)}&& A\to A;\\
&\text{(A2)}&& (A\to B)\to((C\to A)\to(C\to B));\\
&\text{(A3)}&& A\to((A\rightsquigarrow B)\to B);\\
&\text{(A4)}&& (A\rightsquigarrow(B\to C))\to(B\to(A\rightsquigarrow C));\\
&\text{(A5)}&& B\to(A\to A\cdot B);\\
&\text{(A6)}&& (B\to(A\to C))\to(A\cdot B\to C);\\
&\text{(A7)}&& (B\rightsquigarrow B\cdot(B\to A))\to(B\rightsquigarrow A);\\
&\text{(A8)}&& (A\wedge t)\cdot(B\wedge t)\to A\wedge B;\\
&\text{(A9)}&& A\wedge B\to B;\\
&\text{(A10)}&& A\wedge B\to A;\\
&\text{(A11)}&& (C\to A)\wedge(C\to B)\to(C\to A\wedge B);\\
&\text{(A12)}&& A\to A\vee B;\\
&\text{(A13)}&& B\to A\vee B;\\
&\text{(A14)}&& (A\to C)\wedge(B\to C)\to(A\vee B\to C);\\
&\text{(A15)}&& t;\\
&\text{(A16)}&& A\to(t\to A);\\
&\text{(A17)}&& A\to\top;\\
&\text{(A18)}&& \bot\to A;\\
&\text{(PRL)}&&
 \lambda_C((A\vee B)\to A)\vee\rho_D((A\vee B)\to B);\\
&\text{(WCM)}&& (A\rightsquigarrow t)\to(A\to t).
\end{align*}
The four rules, with their premise above their conclusion, are
\[
 \frac{A\qquad A\to B}{B}\;\mathrm{MP},
 \qquad
 \frac{A}{A\wedge t}\;\mathrm{ADJ}_{\mathrm U},
 \qquad
 \frac{A}{B\to A\cdot B}\;\mathrm{PN}_{\to},
 \qquad
 \frac{A}{B\rightsquigarrow B\cdot A}\;\mathrm{PN}_{\rightsquigarrow}.
\]
Proofs are finite trees generated only by these schemata and rules.
An additional designated constant $f$, when present, is interpreted independently
of $t,\bot,\top$ and is fixed by substitution for ordinary variables.

\begin{theorem}[Explicit theorem correspondence]\label{thm:hilbert-explicit}
For every formula $F$,
\[
 \mathrm H_{\mathrm{sl}}^{*}\vdash F
 \quad\Longleftrightarrow\quad \G\vdash\seq F.
\]
Moreover, given a $\G$ proof of any finite hypersequent $H$, there exists
a finite $\mathrm H_{\mathrm{sl}}^{*}$ proof of its negative disjunctive
reading, with repeated components and ordered antecedents retained.
\end{theorem}
\begin{proof}
In the forward direction, construct ordinary proofs of A1--A18 and WCM.
For PRL, communication of identities and the lattice rules give
$(A\vee B\seq A)\mid(A\vee B\seq B)$. Apply left-residual right
introduction to both components, then form the left conjugate with $C$
and the right conjugate with $D$ separately. Join introductions followed
by external contraction give the required singleton conclusion. The four
Hilbert rules translate recursively using the corresponding ordinary rules
and cuts. Thus this direction is a transformation of finite proof trees.

For the converse, take the Lindenbaum quotient by mutual derivability of
implications in this explicitly specified Hilbert calculus. Apply
Lemma~\ref{lem:relative-lindenbaum} with the empty theory to obtain
a bounded residuated lattice with unit $e$ and weak exchange.
The order $[A]\le[B]$ means that there exists a finite Hilbert proof of
$A\to B$; in particular, $e\le[A]$ exactly when $A$ has such a proof.
No total order on this quotient is assumed.

For a sequent $s=(\Gamma\seq A)$ write
\[
 c(s)=(\Pi\Gamma\backslash[A])\wedge e,
 \qquad R(H)=\bigvee_{s\in H}c(s),\qquad R(\varnothing)=\bot.
\]
Here $\Pi\Gamma$ is the ordered product, with empty product $e$; the
syntactic reading is a right-associated join over the actual list with
bottom as its empty tail. For negative $p,z$, independent PRL implies
\[
 p\vee z=e\ \Longrightarrow\
 \lambda_c(p)\vee z=e\quad\hbox{and}\quad\rho_c(p)\vee z=e.
\]
Indeed, use PRL on $p,z$ with the unused parameter equal to $e$ and
$(p\vee z)\backslash p=p$. Interchanging the two branches gives the
right-conjugate statement. Repeating this reasoning permits distinct
left conjugates of the two branches.

Multiplication distributes over joins. If $p\vee z=q\vee w=e$, with
all four values negative, then
\[
 e\le pq\vee z\vee w.
\]
This follows by expanding $(p\vee z)(q\vee w)$ and bounding every mixed
term by the retained background $z\vee w$. Likewise, if $u\vee v=e$
with $u,v\le e$, then
\[
 e\le up\vee vq\vee z\vee w.
\]
For the latter inequality, use $pq\le p,q$ and expand
$(u\vee v)(pq\vee z\vee w)$. These operations never choose a single
component that must be a theorem globally.

We now give the component inequalities for each rule. In the tables,
$p$ is the reading of the sole active premise, or of the \emph{first}
active premise in the displayed order; $q$ is the reading of the second
active premise. Thus $p$ denotes the receiver in $\Cut$ and
$(\backslash L)$, but the producer in $(/L)$. The tables repeat the
active premises to fix this convention. An entry $r$ in the last column
asserts $r\le c(s)$ for the displayed conclusion $s$. All listed weights
are at most $e$. The notation $\Pi\Gamma$ always means the ordered
product of the values of the formulas in $\Gamma$, with empty product
$e$. To keep the weight column short, write $a=[A]$, $b=[B]$,
$g=\Pi\Gamma$ and $d=\Pi\Delta$ in that column.
In particular, for $r\le e$ we have the useful test
\begin{equation}\label{eq:hilbert-component-test}
 r\le c(\Gamma\seq A)
 \quad\Longleftrightarrow\quad (\Pi\Gamma)r\le[A].
\end{equation}

\begin{table}[htbp]
\centering\small
\caption{The four initial cases and four external cases.}
\label{tab:hilbert-initial-external}
\begin{tabularx}{\linewidth}{@{}l>{\raggedright\arraybackslash}X@{}}
\toprule
Initial or external case & Reading identity or inequality\\
\midrule
$A\seq A$ & $c(A\seq A)=e$\\
$\seq t$ & $c(\seq t)=e$\\
$\Gamma,\bot,\Delta\seq A$ & $c(\Gamma,\bot,\Delta\seq A)=e$\\
$\Gamma\seq\top$ & $c(\Gamma\seq\top)=e$\\
Permutation & $R(H)=R(H')$ when $H'$ permutes the component occurrences of $H$\\
$\EW$ & $R(H)\le R(H\mid S)$\\
$\EC$ & $R(H\mid S\mid S)=R(H\mid S)$\\
$\COM$ & $k_1p\le c(s_1)$ and $k_2q\le c(s_2)$, where
$k_1\vee k_2=e$; see~\eqref{eq:hilbert-com-weights}\\
\bottomrule
\end{tabularx}
\end{table}

\begin{table}[htbp]
\centering\small
\caption{The nine unary component cases. The first line of each middle
entry is the premise; the second is the conclusion.}
\label{tab:hilbert-unary}
\begin{tabularx}{\linewidth}{@{}l>{\raggedright\arraybackslash}Xl@{}}
\toprule
Rule & Active premise and conclusion & Lower bound\\
\midrule
$(\cdot L)$ & $\Gamma,A,B,\Delta\seq E$\newline
$\Gamma,A\cdot B,\Delta\seq E$ & $p$\\[3pt]
$(\backslash R)$ & $A,\Sigma\seq B$\newline
$\Sigma\seq A\backslash B$ & $p$\\[3pt]
$(/R)$ & $\Sigma,A\seq B$\newline
$\Sigma\seq B/A$ & $\rho_a(p)$\\[3pt]
$(\wedge L_1)$ & $\Gamma,A,\Delta\seq E$\newline
$\Gamma,A\wedge B,\Delta\seq E$ & $p$\\[3pt]
$(\wedge L_2)$ & $\Gamma,B,\Delta\seq E$\newline
$\Gamma,A\wedge B,\Delta\seq E$ & $p$\\[3pt]
$(\vee R_1)$ & $\Gamma\seq A$\newline
$\Gamma\seq A\vee B$ & $p$\\[3pt]
$(\vee R_2)$ & $\Gamma\seq B$\newline
$\Gamma\seq A\vee B$ & $p$\\[3pt]
$(tL)$ & $\Gamma,\Delta\seq E$\newline
$\Gamma,t,\Delta\seq E$ & $p$\\[3pt]
$\WCM$ & $\Gamma,\Delta\seq t$\newline
$\Delta,\Gamma\seq t$ & $\lambda_g(p)$\\
\bottomrule
\end{tabularx}
\end{table}

\begin{table}[htbp]
\centering\small
\caption{The six binary cases with one active conclusion. External
contexts are retained in their entirety in every row.}
\label{tab:hilbert-binary}
\begin{tabularx}{\linewidth}{@{}l>{\raggedright\arraybackslash}Xl@{}}
\toprule
Rule & First premise; second premise; conclusion & Lower bound\\
\midrule
$(\cdot R)$ & $\Gamma\seq A$; $\Delta\seq B$\newline
$\Gamma,\Delta\seq A\cdot B$ & $\lambda_d(p)q$\\[3pt]
$(\backslash L)$ & $\Gamma,B,\Delta\seq E$; $\Sigma\seq A$\newline
$\Gamma,\Sigma,A\backslash B,\Delta\seq E$
& $\lambda_{(a\backslash b)d}(q)p$\\[3pt]
$(/L)$ & $\Sigma\seq A$; $\Gamma,B,\Delta\seq E$\newline
$\Gamma,B/A,\Sigma,\Delta\seq E$ & $\lambda_d(p)q$\\[3pt]
$(\wedge R)$ & $\Gamma\seq A$; $\Gamma\seq B$\newline
$\Gamma\seq A\wedge B$ & $pq$\\[3pt]
$(\vee L)$ & $\Gamma,A,\Delta\seq E$; $\Gamma,B,\Delta\seq E$\newline
$\Gamma,A\vee B,\Delta\seq E$ & $pq$\\[3pt]
$\Cut$ & $\Gamma,A,\Delta\seq B$; $\Sigma\seq A$\newline
$\Gamma,\Sigma,\Delta\seq B$ & $\lambda_d(q)p$\\
\bottomrule
\end{tabularx}
\end{table}

Here are the ordered calculations underlying the tables. They use
\[
 d\lambda_d(r)\le rd,\qquad \rho_d(r)d\le dr,
 \qquad \lambda_d(r),\rho_d(r)\le e,
\]
which follow directly from their definitions and residuation.
Write $g=\Pi\Gamma$, $d=\Pi\Delta$, $s=\Pi\Sigma$ and
$a=[A]$, $b=[B]$, $h=[E]$ in each applicable case.

For $(\cdot R)$, the premise tests are $gp\le a$ and $dq\le b$.
Consequently
\[
 gd\lambda_d(p)q\le gp\,dq\le ab.
\]
For $(\backslash R)$, the premise test
$asp\le b$ is equivalent to $sp\le a\backslash b$, so the unchanged
weight $p$ works. For $(/R)$, the premise test is $sap\le b$, and
\[
 s\rho_a(p)a\le sap\le b.
\]
Hence $s\rho_a(p)\le b/a$, as required. The apparent asymmetry of the
two right rules comes from placing the weight to the right of the
whole antecedent in~\eqref{eq:hilbert-component-test}.

For $(\backslash L)$, the first premise gives $gbdp\le h$ and the
second gives $sq\le a$. Put $k=\lambda_{(a\backslash b)d}(q)$.
Then the conclusion multiplied by its listed weight satisfies
\[
 \begin{aligned}
 gs(a\backslash b)d\,kp
 &\le gsq(a\backslash b)dp\\
 &\le ga(a\backslash b)dp
 \le gbdp\le h.
 \end{aligned}
\]
For $(/L)$, the first premise instead gives $sp\le a$ and the
second gives $gbdq\le h$. With $k=\lambda_d(p)$ we obtain
\[
 \begin{aligned}
 g(b/a)sd\,kq
 &\le g(b/a)sp\,dq\\
 &\le g(b/a)a\,dq
 \le gbdq\le h.
 \end{aligned}
\]
For $\Cut$, the first premise gives $gadp\le b$ and the second
gives $sq\le a$. With $k=\lambda_d(q)$ the calculation is
\[
 gsd\,kp\le gsqdp\le gadp\le b.
\]
This covers an arbitrary cut formula $A$; WCT is its instance $A=t$.
Every displayed product uses the original prefix, inserted word and
suffix in that order.

For $(\wedge R)$, $pq\le p,q$ implies both $gpq\le a$ and
$gpq\le b$, whence $gpq\le a\wedge b$. For $(\vee L)$, the same
inequalities and join preservation of multiplication give
\[
 g(a\vee b)dpq=gadpq\vee gbdpq\le h.
\]
The unary lattice cases follow by monotonicity; $(\cdot L)$ and
$(tL)$ have identical antecedent values before and after the rule.
The initial cases use reflexivity, the unit law, absorption by the
bottom, and the greatest-element property, respectively.

For WCM the premise gives $gdp\le e$. Unit cyclicity applied to
$g$ and $dp$ gives $dpg\le e$. Thus
\[
 dg\lambda_g(p)\le dpg\le e,
\]
which proves the row in Table~\ref{tab:hilbert-unary} without
exchanging any factors below a general succedent.

For general communication, write the two active premises as
$\Gamma_1,\Xi_1,\Delta_1\seq A$ and
$\Gamma_2,\Xi_2,\Delta_2\seq B$, in that order. Put
$g_i=\Pi\Gamma_i$, $d_i=\Pi\Delta_i$, $x=\Pi\Xi_1$ and
$y=\Pi\Xi_2$, and define
\begin{equation}\label{eq:hilbert-com-weights}
 \begin{gathered}
 u=((x\vee y)\backslash x)\wedge e,\qquad
 v=((x\vee y)\backslash y)\wedge e,\\
 k_1=\lambda_{d_1}(u),\qquad k_2=\lambda_{d_2}(v).
 \end{gathered}
\end{equation}
PRL with both parameters $e$ gives $u\vee v=e$. Moreover,
$yu\le(x\vee y)u\le x$ and $xv\le(x\vee y)v\le y$.
The independent-conjugation closure proved above gives
$k_1\vee k_2=e$. Since $g_1xd_1p\le[A]$ and
$g_2yd_2q\le[B]$, we have
\[
 \begin{aligned}
 g_1yd_1k_1p&\le g_1yu\,d_1p\le g_1xd_1p\le[A],\\
 g_2xd_2k_2q&\le g_2xv\,d_2q\le g_2yd_2q\le[B].
 \end{aligned}
\]
By~\eqref{eq:hilbert-component-test}, these are precisely the two
component bounds in Table~\ref{tab:hilbert-initial-external}.
If $z,w$ are the readings of the two full passive hypersequents,
the earlier choice inequality gives
\[
 e\le k_1p\vee k_2q\vee z\vee w
   \le R\bigl(G_1\mid G_2\mid
       \Gamma_1,\Xi_2,\Delta_1\seq A\mid
       \Gamma_2,\Xi_1,\Delta_2\seq B\bigr).
\]
All these calculations allow empty prefixes, suffixes and exchange
blocks, since their products then equal $e$.

For every unary case, the conjugation closure or the unchanged cover
$p\vee z=e$ transports its weight while retaining $z$. For every
binary single-conclusion case, conjugate the indicated premise cover
first, if necessary, and apply the product-cover inequality in the
listed order. Both passive readings $z,w$ are retained. The separate
choice inequality just checked handles COM. Thus the component bounds
establish the required whole-hypersequent inequality, not merely a
bound for an active component in isolation.

Induction on the actual ordinary proof now proves $e\le R(H)$ for all
23 constructors. Initial rules, permutation, external weakening and
external contraction follow from their lattice readings; the remaining
cases were just established. The proof applies to both the arbitrary-cut
and unit-cut calculi. For $H=(\seq F)$, singleton reflection gives
$e\le[F]$, hence the existence of a finite Hilbert proof of $F$.
\end{proof}

\begin{corollary}[Weak standard completeness for the explicit presentation]
\label{cor:hilbert-weak-standard}
For countably many variables,
$\mathrm H_{\mathrm{sl}}^{*}\vdash F$ if and only if $F$ is valid in
every standard algebra. The same holds with a distinguished, independently
interpreted constant $f$.
\end{corollary}
\begin{proof}
Compose Theorem~\ref{thm:hilbert-explicit} with
Corollary~\ref{cor:fss-weak-standard}. For a countably infinite
alphabet, the bijection $f\mapsto 0$, $p_n\mapsto n+1$ transports
formulas and density-free proofs in both directions. All ordinary rule
instances are preserved under this bijection. A valuation of the
coded atoms corresponds to precisely one choice of the value of $f$
and one valuation of the ordinary variables, as in
Section~\ref{sec:strong-standard}. More generally, for any countable
$V$ the assertion follows from Theorem~\ref{thm:strong-standard} with
$\Gamma=\varnothing$.
\end{proof}

\paragraph{The density rule in the pointed language.}
In a pointed hypersequent calculus a density inference may introduce a
fresh \emph{ordinary variable} $p\in V$. On
$\operatorname{Option}(V)$ its characteristic atom must therefore have
the form $\operatorname{some}(p)$ and meet all the freshness conditions
of~\eqref{eq:main-density}. In the semantic justification one can
then vary $v(p)$ while keeping the interpretation of $f$ fixed.

A generic density calculus over the coded alphabet
$\operatorname{Option}(\mathbb N)$ also permits
$\operatorname{none}$ as the characteristic atom when it is fresh.
It is therefore a different rule system from the pointed calculus.
Nevertheless, their provable hypersequents coincide, as follows. Every proof using only the restricted density
instances is a generic coded proof. Conversely, apply
Corollary~\ref{cor:density-elimination} over the coded alphabet to any
generic coded proof with end-hypersequent $H$. Its density-free proof
of $H$ decodes to an ordinary pointed $\G$ proof, which is also a
proof using only the restricted density instances. Thus the two
density calculi have the same provable hypersequents without additional
hypotheses. The restriction matters for validity in a fixed pointed algebra. For
example, in the standard chain with product $\min$, $e=1$ and $f=0$,
the premise $t\seq f\mid f\seq\bot$ holds under every valuation of
the ordinary variables. Taking $f$ as the characteristic atom would
give the false conclusion $t\seq\bot$.

The arbitrary-theory extension for this same Hilbert presentation is
Theorem~\ref{thm:strong-standard}; the correspondence with $\G$ above
continues to assert formula-theorem equivalence.


\section*{Data and code availability}
This article reports theoretical results. The Lean~4 source code,
dependency specifications, verification scripts, and statement-level
correspondence documentation are available in the archived software
release \emph{GpsUL-1.0.0}, version~1.0.0~\cite{Wang2026GpsUL}.

\section*{Declaration of generative AI and AI-assisted technologies
in the manuscript preparation process}

During the preparation of this manuscript, the author used ChatGPT
(OpenAI) to assist with literature and manuscript organization,
English drafting and revision, and LaTeX preparation.
Its use in mathematical proof development, critical checking,
and the accompanying Lean~4 formalization is described in
Section~\ref{sec:scope}.
After using this tool, the author reviewed and edited the content
as needed and takes full responsibility for the content of the
published article, including its mathematical arguments and references.

\end{document}